%&amstex          
\input amstex\documentstyle{amsppt}  
\pagewidth{12.5cm}\pageheight{19cm}\magnification\magstep1
\topmatter
\title A new basis for the representation ring of a Weyl group\endtitle
\author G. Lusztig\endauthor
\address{Department of Mathematics, M.I.T., Cambridge, MA 02139}\endaddress
\thanks{Supported by NSF grant DMS-1566618.}\endthanks
\endtopmatter   
\document

\define\bx{\boxed}

\define\Irr{\text{\rm Irr}}

\define\mpb{\medpagebreak}

\define\frl{\forall}
\define\pe{\perp}
\define\si{\sim}

\define\sqc{\sqcup}

\define\qua{\quad}

\define\op{\oplus}
   
\define\part{\partial}
\define\emp{\emptyset}

\define\ra{\rangle}
\define\n{\notin}

\define\m{\mapsto}
\define\do{\dots}
\define\la{\langle}

\define\sub{\subset}    

\define\T{\times}
\define\ti{\tilde}
\define\nl{\newline}
\redefine\i{^{-1}}

\define\un{\underline}
\define\ov{\overline}

\define\Hom{\text{\rm Hom}}

\define\tr{\text{\rm tr}}

\redefine\spa{\spadesuit}

\redefine\b{\beta}

\define\g{\gamma}
\redefine\d{\delta}
\define\e{\epsilon}

\define\io{\iota}

\define\ps{\psi}
\define\r{\rho}
\define\s{\sigma}

\define\th{\theta}
\define\k{\kappa}

\define\x{\xi}

\redefine\G{\Gamma}

\define\Th{\Theta}

\define\Ps{\Psi}

\define\BB{\bold B}
\define\CC{\bold C}

\define\FF{\bold F}

\define\NN{\bold N}

\define\QQ{\bold Q}
\define\RR{\bold R}

\define\ZZ{\bold Z}

\define\cc{\Cal C}

\define\ce{\Cal E}
\define\cf{\Cal F}
\define\cg{\Cal G}

\define\ci{\Cal I}

\define\ck{\Cal K}
\define\cl{\Cal L}

\define\car{\Cal R}
\define\cs{\Cal S}

\define\fF{\frak F}

\define\tx{\ti x}

\define\tB{\ti B}

\define\tF{\ti F}

\define\tI{\ti I}

\define\sha{\sharp}

\define\EEIGHT{L1}
\define\CLASSI{L2}
\define\SYMPL{L3}
\define\CLASSII{L4}
\define\CELL{L5}
\define\ORA{L6}
\define\LEA{L7}

\head Introduction and statement of results\endhead
\subhead 0.1\endsubhead   
Let $W$ be an irreducible Weyl group. 
Let $\car_W$ be the (abelian) category of finite dimensional representations of $W$
over $\QQ$ and let $\ck_W$ be the Grothendieck group of $\car_W$. Now $\ck_W$ has a $\ZZ$-basis $\Irr_W$ 
consisting of the irreducible representations of $W$ up to isomorphism. (We often identify a representation of 
$W$ with its isomorphism class.)

Recall that $\Irr_W$ is partitioned into subsets called {\it families}, see
\cite{\EEIGHT, \S8}, \cite{\ORA, 4.2}; these are in 1-1 correspondence with the two-sided cells of $W$. 
For each family $c$ of $W$ we denote by $\car_c$ the (abelian) category of all $E\in\car_W$ 
which are direct sums of irreducible representations in $c$. Let $\ck_c$ the Grothendieck group of $\car_c$.
It has a $\ZZ$-basis consisting of the irreducible representations in $c$. Thus we have 
$\ck_W=\op_c\ck_c$ where $c$ runs over the families of $W$. We now fix a family $c$ of $W$. 

In \cite{\CLASSI} I introduced a class of irreducible objects of $\car_W$ denoted by $\cs_W$ (later called 
special representations); exactly one of these irreducible objects, denoted by $E_c$, is contained in $c$.

In \cite{\CLASSII} I introduced a class of (not necessarily irreducible)
objects of $\car_c$ called ``cells'' (later these objects were called the constructible representations).
In \cite{\CELL} I showed that the constructible representations in $\car_c$ are precisely the representations of 
$W$ carried by the various left cells of $W$ contained in $c$.

In this paper we introduce a class $\BB_c$ of objects of $\car_c$ which includes both $E_c$ and the
constructible representations in $\car_c$ and which forms a $\ZZ$-basis of the group $\ck_c$.
The representations in $\BB_c$ are called {\it new representations}.
(Taking disjoint union over all families of $W$ we obtain a new $\ZZ$-basis of $\ck_W$.)

\subhead 0.2\endsubhead   
Let $\G$ be a finite group. As in \cite{\EEIGHT} we define $M(\G)$ to be the set of all pairs $(x,\r)$
where $x\in\G$ and $\r\in\Irr(Z(x))$ where $Z(x)$ is the centralizer of $x$ in $\G$ and $\Irr(Z(X))$
is the set of irreducible representations of $Z(x)$ over $\CC$ up to isomorphism; these pairs are 
taken up to conjugacy by any element of $\G$.
Let $\CC[M(\G)]$ be the $\CC$-vector space with basis $\{(x,\r);(x,\r)\in M(\G)\}$.

Let $H$ be a subgroup of $\G$. For $x\in\G$ let $(\G/H)^x$ be the 
fixed point set of the left translation action of $x$ on $\G/H$ and let $\CC[(\G/H)^x]$ 
be the $\CC$-vector space with basis $(\G/H)^x$. Now $Z(x)$ acts by
left translation on $(\G/H)^x$ and this induces a linear action of $Z(x)$ on $\CC[(\G/H)^x]$. If 
$\r\in\Irr(Z(x))$, let $N_{H,H,x,\r}$ be the multiplicity of $\r$ in 
the $Z(x)$-module $\CC[(\G/H)^x]$. Let 
$$S_{H,H}=\op_{(x,\r)\in M(\G)}N_{H,H,x,\r}(x,\r)\in\CC[M(\G)].\tag a$$
More generally, let $H\sub H'$ be subgroups of $\G$ with $H$ normal in $H'$. 
Then the obvious surjective map
$\G/H@>>>\G/H'$ restricts to a map $(\G/H)^x@>>>(\G/H')^x$ and this induces a 
linear map $\CC[(\G/H)^x]@>>>\CC[(\G/H')^x]$ (compatible with $Z(x)$ actions) whose image 
is denoted by $\ci$. Now $\ci$ is a $Z(x)$-submodule of $\CC[(\G/H')^x]$. If $\r\in\Irr(Z(x))$,
let $N_{H,H',x,\r}$ be the multiplicity of $\r$ in the $Z(x)$-module $\ci$. Let 
$$S_{H,H'}=\op_{(x,\r)\in M(\G)}N_{H,H',x,\r}(x,\r)\in\CC[M(\G)].\tag b$$
For example,
$$S_{\{1\},\{1\}}=\sum_{\r\in\Irr(\G)}\dim\r(1,\r),$$
$$S_{\{1\},\G}=(1,1),$$
$$S_{\G,\G}=\sum_{x\in\G\text{ up to conjugacy}}(x,1).$$

\subhead 0.3\endsubhead   
As in \cite{\ORA, \S4} we attach to $c$ a finite group $\cg_c$ and an imbedding $c@>>>M(\cg_c)$.
Let $M_0(\cg_c)$ be the image of this imbedding. For $(x,\r)\in M_0(\cg_c)$ let $E_{x,\r}$ be the corresponding 
(irreducible) representation in $c$. For any $\ce\in\car_c$ we define $\un\ce\in\CC[M(\cg_c)]$ by 
$\un\ce=\sum_{(x,\r)\in M_0(\cg_c)}(E_{x,\r}:\ce)(x,\r)$
where $(E_{x,\r}:\ce)\in\NN$ is the multiplicity of $E_{x,\r}$ in $\ce$. Note that $\ce\m\un\ce$ defined an 
imbedding $\ck_c@>>>\CC[M(\cg_c)]$.

As was pointed out in \cite{\LEA}, to any constructible representation $E$ in $\car_c$ one can attach a subgroup 
$H_E$ of $\cg_c$, well defined up to conjugacy, such that $\un E=S_{H_E,H_E}$, see 0.2(a). Moreover, 

(a) $E\m H_E$ 
\nl
is an injective map from the set of constructible representations in $\car_c$ to the set of 
subgroups of 
$\cg_c$ (up to conjugacy). Let $\fF_c$ be the set of subgroups of $\cg_c$ which are conjugate to a subgroup in 
the image of the map (a). We have $\cg_c\in\fF_c$.
We say that $c$ is {\it anomalous} if $\{1\}\n\fF_c$. 
If $W$ is of classical type then $c$ is not anomalous. If $W$ is of exceptional type, then $c$ 
is anomalous in exactly the following cases:

(b) the unique $c$ with $|c|=2$ with $W$ of type $E_7$;

(c) the two $c$ with $|c|=2$ with $W$ of type $E_8$;  

(d) the unique $c$ with $|c|=4$ with $W$ of type $G_2$;

(e) the unique $c$ with $|c|=11$ with $W$ of type $F_4$;

(f) the unique $c$ with $|c|=17$ with $W$ of type $E_8$.
\nl
Let $\hat\fF_c$ be the set of subgroups of $\cg_c$ which are either $\{1\}$ or are in $\fF_c$.
Let $\ti\Th_c$ be the set of all pairs $(H,H')$ where $H\in\hat\fF_c,H'\in\hat\fF_c$ and $H$ is a 
normal subgroup of $H'$. Now $\cg_c$ acts on $\ti\Th_c$ by simultaneous conjugation. 
We now state our main result.

\proclaim{Theorem 0.4} There exists a $\cg_c$-stable subset $\Th_c$ of $\ti\Th_c$ such that the following hold:

(i) For any $H\in\fF_c$ we have $(H,H)\in\Th_c$.

(ii) We have $(1,\cg_c)\in\Th_c$.

(iii) For any $(H,H')\in\Th_c$ there is a unique object $E_{H,H'}\in\car_c$ such that
$S_{H,H'}=\un E_{H,H'}$, see 0.2(a). Let $\BB_c$ be the set of isomorphism classes of objects 
of $\car_c$ of the form $E_{H,H'}$ for some $(H,H')\in\Th_c$.

(iv) The map $(H,H')\m E_{H,H'}$ defines a bijection from the set of $\cg_c$-orbits on $\Th_c$
to $\BB_c$. Moreover $\BB_c$ is a $\ZZ$-basis of $\ck_c$.
\endproclaim
The representations in $\BB_c$ are the new representations mentioned in 0.1.
From (i) we see that any constructible representation of $\car_c$ is in $\BB_c$.
From (ii) we see that the special representation $E_c$ is in $\BB_c$.

In the case where $W$ is of type $A$ the theorem is trivial; we have $\cg_c=\{1\}$ and $\BB_c$ 
consists of the unique representation in $c$. The proof of the theorem for $W$ of type $B_n,C_n,D_n$ 
is given in \S2.
The proof of the theorem for $W$ of exceptional type is given in \S3.

\subhead 0.5\endsubhead
In this paper we also define a canonical bijection $c@>\si>>\BB_c$, $E\m\hat E$ which has the property 
that for any $E\in c$, $E$ has appears with multiplicity one in $\hat E$.
For $E,E'$ in $c$ let $E':\hat E$ be the multiplicity of $E'$ in $\hat E$.
Property (i) below will be proved in a sequel to this paper. (For $W$ of exceptional type (i) is easily 
deduced from the formulas in 3.2-3.8.)

(i) The matrix $(E':\hat E)$ indexed by $c\T c$ is upper triangular unipotent for a suitable partial order on $c$.

\subhead 0.6\endsubhead
In the setup of 0.2 we define (following \cite{\EEIGHT, \S4}) a pairing $\{,\}:M(\G)\T M(\G)@>>>\CC$ by 
$$\align&\{(x,\r),(x',\r')\}\\&
=|Z(x)|\i|Z(x')|\i\sum_{g\in\G; x gx'g\i=gx'g\i x}\ov{\tr(g\i xg,\r')}\tr(gx'g\i,\r)\endalign$$
where $\bar{}$ is complex conjugation.
We define the non-abelian Fourier transform $A:\CC[M(\G)]@>>>\CC[M(\G)]$ as the $\CC$-linear map
such that $$A(x,\r)=\sum_{(x',\r')\in M(\G)}\{(x,\r),(x',\r')\}(x',\r')$$ for any $(x,\r)\in M(\G)$. 
According to \cite{\EEIGHT}, we have $A^2=1$. Let $M(\G)_{\ge0}$ be the set of elements
$$\sum_{(x,\r)\in M(\G)}c_{x,\r}(x,\r)\in\CC[M(\G)]$$ such that $c_{x,\r}\in\RR_{\ge0}$ for any $(x,\r)\in M(\G)$.

An element $f\in\CC[M(\G)]$ is said to be {\it bipositive} if $f\in M(\G)_{\ge0}$ and $A(f)\in M(\G)_{\ge0}$.
We have the following result.

\proclaim{Theorem 0.7} Let $H\sub H'$ be subgroups of $\G$ with $H$ normal in $H'$. Then
$S_{H,H'}\in\CC[M(\G)]$ is bipositive. Hence (by 0.4), if $\G=\cg_c$ and $\ce$ is a new representation 
in $\car_c$, then $\un\ce\in\CC[M(\G)]$ is bipositive.
\endproclaim
The proof is given in \S4.

\subhead 0.8\endsubhead
In a sequel to this paper we will extend the results of the paper by constructing a new basis for 
$\CC[M(\cg_c)]$ consisting of bipositive elements; this provides a new $\ZZ$-basis for the
Grothendieck group of unipotent representations of a split Chevalley group over a finite field.

\subhead 0.9\endsubhead
{\it Notation.} For $a\le b$ in $\NN$ we write $[a,b]=\{z\in\NN;a\le z\le b\}$. We set $[1,0]=\emp$. 
For a finite set $Y$ we write $|Y|$ for the cardinal of $Y$. For $a,b$ in $\ZZ$ we write $a=_2b$ if
$a=b\mod2$ and $a\ne_2b$ if $a\ne b\mod2$. We write $\ZZ/2\ZZ=\FF_2$.

\head 1. The set $S_D$\endhead  
\subhead 1.1\endsubhead
Let $D\in\NN$. A subset $I$ of $[1,D]$ is said to be an {\it interval} if $I=[a,b]$ for some 
$a\le b$ in $[1,D]$. Let $\ci_D$ be the set of intervals of 
$[1,D]$. For $I=[a,b],I'=[a',b']$ in $\ci_D$ we write $I\prec I'$ whenever $a'<a\le b<b'$. 
We say that $I,I'$ are non-touching (and we write $I\spa I'$) if $a'-b\ge2$ or $a-b'\ge2$.
Let $\ci_D^1=\{I\in\ci_D;|I|=\text{ odd}\}$.
Let $R^1_D$ be the set whose elements are the subsets of $\ci^1_D$. 
Let $\emp\in R^1_D$ be the empty subset of $\ci^1_D$.

When $D\ge2$ and $i\in[1,D]$ we define an (injective) map $\x_i:\ci_{D-2}@>>>\ci_D$ as follows:
$$\align&\x_i([a',b'])=[a'+2,b'+2]\text{ if }i\le a',\qua
\x_i([a',b'])=[a',b']\text{ if }i\ge b'+2, \\&
\x_i([a',b'])=[a',b'+2]\text{ if }a'<i<b'+2.\tag a\endalign$$ 
We have $\x_i(\ci^1_{D-2})\sub\ci^1_D$. 
We define $t_i:R^1_{D-2}@>>>R^1_D$ 
by $B'\m\{\x_i(I');I'\in B'\}\sqc\{i\}$. We have $|t_i(B')|=|B'|+1$.

\subhead 1.2\endsubhead
We define a subset $S_D$ of $R^1_D$ by induction on $D$ as follows. When $D\in\{0,1\}$, $S_D$ 
consists of a single element, namely $\emp\in R^1_D$. When $D\ge2$ we say that $B\in R^1_D$ is in 
$S_D$ if either $B=\emp$ or if

(i) there exists $i\in[1,D]$ (if $D$ is even) or $i\in[1,D-1]$ (if $D$ is odd) and $B'\in S_{D-2}$ such that 
$B=t_i(B')$. 
\nl
If $D$ is odd, we have $S_D=S_{D-1}$ (use induction on $D$).

Until the end of 1.8 we assume that $D$ is even.

\subhead 1.3. The set $S'_D$\endsubhead
Let $B\in R^1_D$. We consider the following properties $(P_0),(P_1)$ that $B$ may or may not have.

$(P_0)$ {\it If $I\in B$, $\tI\in B$, then either $I=\tI$, or $I\spa\tI$, or $I\prec\tI$, or 
$\tI\prec I$.}

$(P_1)$ {\it If $[a,b]\in B$ and $c\in\NN$ satisfies $a<c<b$, $a-c=_21$ (hence $b-c=_21$), then 
there exists $[a_1,b_1]\in B$ such that $a<a_1\le c\le b_1<b$.}
\nl
From the definitions we see that if $D\ge2$, $i\in[1,D]$, $B'\in R^1_{D-2}$ and $B=t_i(B')\in R^1_D$, the
following holds.

(a) {\it $B'$ satisfies $(P_0)$ if and only if $B$ satisfies $(P_0)$; $B'$ satisfies $(P_1)$ if and 
only if $B$ satisfies $(P_1)$.}
\nl
Let $S'_D$ be the set of all $B\in R^1_D$ which satisfy $(P_0),(P_1)$. 
In the setup of (a) we have the following consequence of (a).

(b) {\it We have $B'\in S'_{D-2}$ if and only if $B\in S'_D$.}
\nl
We show:

(c) {\it $S_D=S'_D$.}
\nl
We argue by induction on $D$. If $D=0$, $S'_D$ consists of the empty set hence (c) holds in this case. 
Assume now that $D\ge2$. Let $B\in S_D$. We show that $B\in S'_D$. If $B=\emp$ this is clear.
If $B\ne\emp$ then $B=t_i(B')$ for some 
$i,B'\in S_{D-2}$. By the induction hypothesis we have $B'\in S'_{D-2}$. By (b) we have
$B\in S'_D$. We see that $B\in S_D\implies B\in S'_D$. Conversely, let $B\in S'_D$. We show that 
$B\in S_D$. If $B=\emp$ this is obvious. Thus we can assume that $B\ne\emp$. 
Let $[a,b]\in B$ be such that $b-a$ is minimum. If $a<z<b$, 
$z=_2a+1$ then by $(P_1)$ we have $z\in[a',b']$ with $[a',b']\in B$, $b'-a'<b-a$, contradicting 
the minimality of $b-a$. We see that no $z$ as above exists. Thus, $[a,b]=\{i\}$ for some $i\in[1,D]$. 
Using $(P_0)$ and $\{i\}\in B$, we see that $B$ does not contain any interval of the form $[a,i]$ with 
$a<i$, or $[i,b]$ with $i<b$, or $[a,i-1]$ with $a<i$ or $[i+1,b]$ with $i<b$; hence any interval of $B$
other than $\{i\}$ is of the form $\x_i[a',b']$ where $[a',b']\in\ci^1_{D-2}$. Thus we have $B=t_i(B')$ 
for some $B'\in S_{D-2}$. From (b) we see that $B'\in S'_{D-2}$. Using the induction 
hypothesis we deduce that $B'\in S_{D-2}$. By the definition of $S_D$, we have $B\in S_D$. This 
completes the proof of (c).

\mpb

The following result has already been proved as a part of the proof of (c).

(d) {\it Assume that $D\ge2$, $i\in[1,D]$. Let $B\in S_D$ be such that $\{i\}\in B$. Then there exists
$B'\in S_{D-2}$ such that $B=t_i(B')$.}

\subhead 1.4\endsubhead
For $B\in S_D$, $j\in[1,D]$ we set $B_j=\{I\in B;j\in I\}$. From the definitions we deduce:

(a) {\it Assume that $D\ge2$, $i\in[1,D]$ and that $B'\in S_{D-2},B=t_i(B')\in S_D$. Then for $r\in[1,D-2]$ we 
have:

$|B'_r)|=|B_r|$ if $r\le i-2$, $|B'_r|=|B_{r+2}|$ if $r\ge i$,

$|B_{i-1}|=|B_{i+1}|=|B'_{i-1}|$, $|B_i|=|B'_{i-1}|+1$ if $1<i<D$,

$|B_{i-1}|=0$ if $i=D$, $|B_{i+1}|=0$ if $i=1$.}

\subhead 1.5\endsubhead
Let $B\in S_D, B\ne\emp$. In this case we must have $\{j\}\in B$ for some $j\in[1,D]$; we assume that $j$ is as 
small as possible (then it is uniquely determined). As in the proof of 1.3(c) we have $B=t_j(B')$ where 
$B'\in S_{D-2}$. Let $i$ be the smallest number in $\cup_{I\in B}I$. We have $i\le j$. We show:

(a) {\it For any $h\in[i,j]$, we have $[h,\ti h]\in B$ for a unique $\ti h\in[h,D]$; moreover we have 
$j\le\ti h$.}
\nl
We argue by induction on $D$. When $D=0$ the result is obvious. We now assume that $D\ge2$.
Assume first that $i=j$. By $(P_0)$, $\{j\}\in B$ implies that we cannot have $[j,b]\in B$ with
 $j<b$; thus (a) holds in this case. In particular, (a) holds when $D=2$ (in this case we have $i=j$). We 
now assume that $D\ge4$. We can assume that $i<j$. We have $[i,b]\in B$ for some $b>i$ hence 
$|B|\ge2$ so that $|B'|\ge1$ and $B'\ne\emp$. Then $i',j'$ are defined in terms of 
$B'$ in the same way as $i,j$ are defined in terms of $B$. From $(P_1)$ we see that there exists $j_1$ 
such that $i<j_1<b$ such that $\{j_1\}\in B$. By the minimality of $j$ we must have $j\le j_1$. Thus we 
have $i<j<b$. We have $[i,b]=\x_j[i,b-2]]$ hence $[i,b-2]\in B'$. This implies that $i'\le i$. We 
have $[i',c]\in B'$ for some $c\in[i',D-2]$, $c=_2i'$; hence $[i',c']\in B$ for some $c'\ge i'$ 
so that $i'\ge i$. Thus we have $i'=i$. By the induction hypothesis, the following holds:

(b) {\it For any $r\in[i,j']$, we have $[r,r_1]\in B'$ for a unique $r_1$; moreover $j'\le r_1$.}
\nl
If $j'\le j-2$ then $\{j'\}=\x_j(\{j'\})\in B$. Hence $j'\ge j$ by the minimality of $j$; this is a 
contradiction. Thus we have $j'\ge j-1$.

Let $r\in[i,j-1]$. Then we have also $r\in[i,j']$ hence $r_1$ is defined as in (b). We have
$[r,r_1]\in B'$ hence $[r,r_1+2]\in B$ (we use that $r<j\le j'+1\le r_1+1<r_1+2$); we have 
$j<r_1+2$. Assume now that $[r,r_2]\in B$ with $r\le r_2$. Then $r<r_2$ (by the minimality of $j$).
If $j=r_2$ or $j=r_2+1$ then applying $(P_0)$ to $\{j\},[r,r_2]$ gives a contradiction.
Thus we must have either $r<j<r_2$ or $j>r_2+1$. If $j>r_2+1$ then $[r,r_2]\in B'$ hence by (b), 
$r_2=r_1$, hence $j>r_1+1$ contradicting $j<r_1+2$. Thus we have $r<j<r_2$, so that 
$[r,r_2-2]\in B'$ hence by (b), $r_2-2=r_1$. Thus we have $r<j<r_2$ so that $[r,r_2-2]\in B'$ 
hence by (b), $r_2-2=r_1$.

Next we assume that $r=j$. In this case we have $\{r\}\in B$. Moreover, if $[r,r']\in B$ with
$r\le r'\le D$, then we cannot have $r<r'$ (if $r<r'$ then applying $(P_0)$ to $\{r\},[r,r']$ gives a
contradiction). This proves (a).

We show:

(c) {\it Assume that $j<D$ and that $i\le h<j$. Then $\ti h$ in (a) satisfies $\ti h>j$.}
\nl
Assume that $\ti h=j$, so that $[h,j]\in B$. Since $h<j$, 
applying $(P_0)$ to $\{j\},[h,j]$ gives a contradiction.  This proves (c).

(d) {\it Assume that $j<D$ and that $r\in[j+1,D]$. We have $[j+1,r]\n B$.}
\nl
Assume that $[j+1,r]\in B$. Applying $(P_0)$ to $\{j\},[j+1,r]$ gives a contradiction.  This proves (d).

We show:

(e) {\it For $h\in[i,j]$ we have $|B_h|=h-i+1$. If $j<D$ we have $|B_{j+1}|=j-i$.}
\nl
Let $h\in[i,j]$. Then for any $h'\in[i,h]$, $B_h$ contains $[h',\ti h']$ (since $h\le\ti h'$), see (a).
Conversely, assume that $[a,b]\in B_h$. We have $a\le h$. By the definition of $i$ we have $i\le a$.
By the uniqueness statement in (a) we have $b=\ti a$ so that $[a,b]$ is one of the $h-i+1$ intervals 
$[h',\ti h']$ above. This proves the first assertion of (e). Assume now that $j<D$. If $h'\in[i,j]$, 
$h'<j$,  then $[h',\ti h']\in B_{j+1}$, by (c). Conversely, assume that
$[a,b]\in B_{j+1}$. We have $a\le j+1$ and by (d) we have $a\ne j+1$ so that $a\le j$. If $a=j$, then 
by the uniqueness in (a) we have $b=j$ which contradicts $j+1\in[a,b]$. Thus we have $a\le j-1$. We see 
that $[a,b]$ is one of the $j-i$ intervals $[h',\ti h']$ with $h'\in[i,j]$, $h'<j$. This proves (e).

\subhead 1.6\endsubhead  
For $B\in S_D$, $j\in[1,D]$, we set 
$$\e_j(B)=|B_j|(|B_j|+1)/2\in \FF_2.$$

We have $\e_j(B)=1$ if $|B_j|\in(4\ZZ+1)\cup(4\ZZ+2)$, $\e_j(B)=0$ if $|B_j|\in(4\ZZ+3)\cup(4\ZZ)$. 

Assume now that $B\ne\emp$. Let $i\le j$ in $[1,D]$ be as in 1.5. From 1.5(e) we deduce:

(a) {\it We have $(|B_i|,|B_{i+1}|,\do,|B_j|)=(1,2,3,\do,j-i,j-i+1)$. If 
$j<D$, we have $|B_{j+1}|=j-i$.}
\nl
From (a) we deduce:
$$\align&(\e_i(B),e_{i+1}(B),\do,\e_j(B))\\&=
(1\T2)/2,(2\T3)/2,(3\T4)/2,\do,(j-i)(j-i+1)/2,(j-i+1)(j-i+2)/2);\tag b\endalign$$
(c) {\it if $j<D$, then} $\e_{j+1}(B)=(j-i)(j-i+1)/2$.
\nl
For future reference we note:

(d) {\it If $c\in\ZZ$ then $c(c+1)/2\ne_2(c+2)(c+3)/2$.}

(e) {\it If $c\in2\ZZ$ then $c(c+1)/2\ne_2(c+1)(c+2)/2$.}

\subhead 1.7\endsubhead
Let $B\in S_D$, $\tB\in S_D$ be such that $B\ne\emp,\tB\ne\emp$ and $\e_h(B)=\e_h(\tB)$ for any $h\in[1,D]$. We 
show:

(a) {\it We can find $z\in[1,D]$ such that $\{z\}\in B$, $\{z\}\in\tB$.}
\nl
We associate $i\le j$ to $B$ as in 1.5; let $\ti i\le\ti j$ be the analogous number for $\tB$.
Assume first that $j<\ti j$ (so that $j<D$) and $i<\ti i$. From 1.6 for $B$ we have 
$\e_i(B)=(1\T2)/2=1$. Since $i<\ti i$ we have $\e_i(\tB)=0$. Hence $1=_20$, a contradiction.
Thus we must have $i\ge\ti i$. 

Next we asssume that $j<\ti j$ (so that $j<D$) and $\ti i<i$. From 1.6 for $\tB$ we have 
$\e_{\ti i}(\tB)=(1\T2)/2$; moreover $\e_{\ti i}(B)=0$. Hence $1=_20$, a contradiction.  
Thus when $j<\ti j$ we must have $i=\ti i$. 
From 1.6(c) for $B$ we have $e_{j+1}(B)=(j-i)(j-i+1)/2$  and from 1.6(b) for $\tB$ we 
have $e_{j+1}(\tB)=(j-i+2)(j-i+3)/2$. It follows that 
$$(j-i)(j-i+1)/2)=_2(j-i+2)(j-i+3)/2,$$ 
contradicting 1.6(d). We see that $j<\ti j$ leads to a contradiction. Similarly, $\ti j<j$ leads to a 
contradiction. Thus we must have $j=\ti j$, so that (a) holds with $z=j=\ti j$. This completes the 
proof of (a).

\subhead 1.8\endsubhead
Let $B\in S_D$, $\tB\in S_D$.

(a) {\it Assume that $\tB=\emp$ and that $\e_h(B)=\e_h(\tB)$ for any $h\in[1,D]$. Then $\tB=B$.}
\nl
The proof is similar to that of 1.7(a). Assume that $B\ne\emp$. Let $i\le j$ be attached 
to $B$ as in 1.5. 

Using 1.6 we see that $e_i(B)=(1\T2)/2$. On the other hand we have $e_i(\tB)=0$. We get $1=_20$, a
 contradiction. This proves (a).

\subhead 1.9\endsubhead
We no longer assume that $D$ is even. Let $V$ be the $\FF_2$-vector space consisting of all functions 
$[1,D]@>>>\FF_2$. For any subset $I$ of $[1,D]$ let $e_I\in V$ be the function 
whose value at $i$ is $1$ if $i\in I$ and is $0$ if $i\n I$. For $i\in[1,D]$ we set $e_i=e_{\{i\}}$. Now 
$\{e_i;i\in[1,D]\}$ is a basis of $V$. We define a symplectic form $(,):V\T V@>>>\FF_2$ by 
$(e_i,e_j)=1$ if $i-j=\pm1$, $(e_i,e_j)=0$ if $i-j\ne\pm1$. This symplectic form is nondegenerate if $D$ is even
while if $D$ is odd it has a one dimensional radical spanned by $e_1+e_3+e_5+\do+e_D$.

For any subset $Z$ of $V$ we set $Z^\pe=\{x\in V;(x,z)=0\qua\frl z\in Z\}$.

When $D\ge2$ we denote by $V'$ the $\FF_2$-vector space consisting of all functions $[1,D-2]@>>>\FF_2$.
For any $I'\sub[1,D-2]$ let $e'_{I'}\in V'$ be the function whose value at $i$ is $1$ if $i\in I'$ and 
is $0$ if $i\n I'$. For $i\in[1,D-2]$ we set $e'_i=e'_{\{i\}}$. Now $\{e'_i;i\in[1,D-2]\}$ is a basis of 
$V'$. We define a symplectic form $(,)':V'\T V'@>>>\FF_2$ by $(e'_i,e'_j)=1$ if $i-j=\pm1$,
$(e'_i,e'_j)=0$ if $i-j\ne\pm1$. 

When $D\ge2$, for any $i\in[1,D]$ there is a unique linear map $T_i:V'@>>>V$ such that 
the sequence $T_i(e'_1),T_i(e'_2),\do,T_i(e'_{D-2})$ is:

$e_1,e_2,\do,e_{i-2},e_{i-1}+e_i+e_{i+1},e_{i+2},e_{i+3},\do,e_D$ (if $1<i<D$),

$e_3,e_4,\do,e_D$ (if $i=1$),

$e_1,e_2,\do,e_{D-2}$ (if $i=D$).
\nl
Note that $T_i$ is injective and $(x,y)'=(T_i(x),T_i(y))$ for any $x,y$ in $V'$.
For any $I'\in\ci^1_{D-2}$ we have $T_i(e'_{I'})=e_{\x_i(I')}$. Let $V_i$ be the image of $T_i:V'@>>>V$. 
From the definitions we deduce:

(a) {\it We have $e_i^\pe=V_i\op \FF_2e_i$.}
\nl
We now assume that $D$ is even.
For $j\in[1,D-2]$ let $\e'_j:S_{D-2}@>>>\FF_2$ be the 
analogue of $\e_i:S_D@>>>\FF_2$ when $D$ is replaced by $D-2$.

For $B\in S_D$, we define $\e(B)\in V$ by $i\m\e_i(B)$. 
For $B'\in S_{D-2}$ we define $\e'(B')\in V'$ by $j\m\e'_j(B')$. We show:

(b) {\it Assume that $D\ge2$, $i\in[1,D]$. Let $B'\in S_{D-2},B=t_i(B')\in S_D$. Then
$\e(B)=T_i(\e'(B'))+ce_i$ for some $c\in\FF_2$.}
\nl
An equivalent statement is: for any $j\in[1,D]-\{i\}$ we have $\e_j(B)=\e'_{j'}(B')$ 
if $j'\in[1,D-2]$ is such that $j\in\x_i(\{j'\})$; and $\e_j(B)=0$ if no such $j'$ exists. It is enough to show:

$|B'_h|=|B_h|$ if $h\in[1,i-2]$,

$|B'_{h-2}|=|B_h|$ if $h\in[i+2,D]$,

$|B_{i-1}|=|B_{i+1}|=|B'_{i-1}|$ if $1<i<D$,

$|B_{i-1}|\in\{0,-1\}$ (hence $\e_{i-1}(B)=0$) if $i=D$, 

$|B_{i+1}|\in\{0,-1\}$ (hence $\e_{i+1}(B)=0$) if $i=1$.
\nl
This follows from 1.4(a).

\mpb

For $B\in S_D$ let $\la B\ra$ be the subspace of $V$ generated by $\{e_I;I\in B\}$. 
For $B'\in S_{D-2}$ let $\la B'\ra$ be the subspace of $V'$ generated by $\{e'_{I'};I'\in B'\}$. We show:

(c) {\it Let $B\in S_D$. We have $\e(B)\in\la B\ra$.
If $D\ge2,i\in[1,D]$, $B'\in S_{D-2},B=t_i(B')\in S_D$, then $\la B\ra=T_i(\la B'\ra)\op\FF_2e_i$.}
\nl
To prove the first assertion of (c) we argue by induction on $D$. For 
$D=0$ there is nothing to prove. Assume that $D\ge2$. Let $i,B'$ be as 
in (b). By the induction hypothesis we have $\e'(B')\in\la B'\ra\sub V'$. Using (b) we see that it is 
enough to show that $T_i(\la B'\ra)\sub\la B\ra$. (Since $\{i\}\in B$, we have $e_i\in\la B\ra$.) 
Using the equality $T_i(e'_{I'})=e_{\x_i(I')}$ for any $I'\in B'$ it remains to note that $\x_i(I')\in B$ 
for $I'\in B'$. This proves the first assertion of (c). The same proof shows the second assertion of (c).

\subhead 1.10\endsubhead
Let $B\in S_D,\tB\in S_D$. We show:

(a) {\it If $\e(B)=\e(\tB)$, then $B=\tB$. }
\nl
Wec argue by induction on $D$. If $D=0$, there is nothing to prove. Assume that $D\ge2$. If $\tB=\emp$,
(a) follows from 1.8(a). Similarly, (a) holds if and $B=\emp$. Thus, we can 
assume that $B\ne\emp$, $\tB\ne\emp$. By 1.7(a) we can find $i\in[1,D]$ such that 
$\{i\}\in B$, $\{i\}\in\tB$. By 1.3(d) we then have $B=t_i(B')$, $\tB=t_i(\tB')$ with 
$B'\in S_{D-2}$, $\tB'\in S_{D-2}$. Using our assumption and 1.9(b) we see that 
$T_i(\e'(B'))=T_i(\e'(\tB'))+ce_i$ for some $c\in\FF_2$. Using 1.9(a) we see that $c=0$ so that 
$T_i(\e'(B'))=T_i(\e'(\tB'))$. Since $T_i$ is injective, we deduce $\e'(B')=\e'(\tB')$.
By the induction hypothesis we have $B'=\tB'$ hence $B=\tB$. This proves (a).

\subhead 1.11\endsubhead
Any $x\in V$ can be written uniquely in the form
$$x=e_{[a_1,b_1]}+e_{[a_2,b_2]}+\do+e_{[a_r,b_r]}$$
where $[a_r,b_r]\in\ci_D$ are such that any two of them are non-touching and
$r\ge0$, $1\le a_1\le b_1<a_1\le b_2<\do<a_r\le b_r\le D$.
Following \cite{\SYMPL, 3.3} we set
$$u(v)=|\{s\in[1,r];a_s=_20,b_s=_21\}|-|\{s\in[1,r];a_s=_21,b_s=_20\}|\in\ZZ.\tag a$$
This defines a function $u:V@>>>\ZZ$. When $D\ge2$ we denote by $u':V'@>>>\ZZ$ the analogous
function with $D$ replaced by $D-2$. We show:

(b) {\it Assume that $D\ge2,i\in[1,D]$. Let $v'\in V'$ and let $v=T_i(v')+ce_i\in V$ where $c\in\FF_2$.
We have $u(v)=u'(v')$.}
\nl
We write $v'=e'_{[a'_1,b'_1]}+e'_{[a'_2,b'_2]}+\do+e'_{[a'_r,b'_r]}$ where $r\ge0$, 
$[a'_s,b'_s]\in\ci_{D-2}$ for all $s$ and any two of $[a'_s,b'_s]$ are non-touching. For each $s$, we 
have $T_i(e'_{[a'_s,b'_s]})=e_{[a_s,b_s]}$ where $[a_s,b_s]=\x_i[a'_s,b'_s]$ so that $a_s=_2a'_s$,
$b_s=_2b'_s$ and the various $[a_s,b_s]$ which appear are still non-touching with each other. Hence 
$u(T_i(v'))=u'(v')$. We have $v=T_i(v')$ or $v=T_i(v')+e_i$. If $v=T_i(v')$, we have $u(v)=u'(v')$, as desired. 
Assume now that $v=T_i(v')+e_i$. From the definition of $\x_i$ we see that either

(i) $[i,i]$ is non-touching with any $[a_s,b_s]$, or

(ii) $[i,i]$ is not non-touching with some $[a,b]=[a_s,b_s]$ which is uniquely
determined and we have $a<i<b$.

If (i) holds then $e_i$ does not contribute to $u(v)$ and $u(v)=u(T_i(v'))=u'(v')$. 
We now assume that (ii) holds. Then $e_{[a,b]}+e_i=e_{[a,i-1]}+e_{[i+1,b]}$. We consider six cases.

(1) $a$ is even $b$ is odd, $i$ is even; then $|[i+1,b]|$ is odd so that the contribution of  
$e_{[a,i-1]}+e_{[i+1,b]}$ to $u(v)$ is $1+0$; this equals the contribution of $e_{[a,b]}$ to $u(T_i(v'))$ 
which is $1$.

(2) $a$ is even, $b$ is odd, $i$ is odd; then $|[a,i-1]|$ is odd so that the contribution of  
$e_{[a,i-1]}+e_{[i+1,b]}$ to $u(v)$ is $0+1$; this equals the contribution of $e_{[a,b]}$ to $u(T_i(v'))$ 
which is $1$.

(3) $a$ is odd, $b$ is even, $i$ is even; then $|[i+1,b]|$ is odd so that the contribution of  
$e_{[a,i-1]}+e_{[i+1,b]}$ to $u(v)$ is $0-1$; this equals the contribution of $e_{[a,b]}$ to $u(T_i(v'))$ 
which is $-1$.

(4) $a$ is odd, $b$ is even, $i$ is odd; then $|[a,i-1]|$ is odd so that the contribution of  
$e_{[a,i-1]}+e_{[i+1,b]}$ to $u(v)$ is $-1+0$; this equals the contribution of $e_{[a,b]}$ to $u(T_i(v'))$
 which is $-1$.

(5) $a=_2b=_2i+1$; then $|[a,i-1]|$ is odd, $|[i+1,b]|$ is odd so that the contribution of  
$e_{[a,i-1]}+e_{[i+1,b]}$ to $u(v)$ is $0+0$; this equals the contribution of $e_{[a,b]}$ to $u(T_i(v'))$ 
which is $0$.

(6) $a=_2b=_2i$; then the contribution of $e_{[a,i-1]}+e_{[i+1,b]}$ to $u(v)$ is $1-1$ or $-1+1$; this 
equals the contribution of $e_{[a,b]}$ to $u(T_i(v'))$ which is $0$.
\nl
This proves (b).

\subhead 1.12\endsubhead
We view $V$ as the set of vertices of a graph in which $x,x'$ in $V$ are joined whenever 
there exists $i\in[1,D]$ such that $x+x'=e_i$, $(x,e_i)=(x',e_i)=0$. 
Similarly if $D\ge2$, we view $V'$ as the set of vertices of a graph in which $y,y'$ in $V'$ are joined whenever 
there exists $i\in[1,D-2]$ such that $y+y'=e'_i$, $(y,e'_i)'=(y',e'_i)'=0$. We show: 

(a) {\it If $y,y'$ in $V'$ are joined in the graph $V'$ 
 then $T_i(y),T_i(y')$ are in the same connected component of the graph $V$.}
\nl
We can find $j\in[1,2d-2]$ such that $(y,e'_j)'=(y',e'_j)'=0$, $y+y'=e'_j$.
Hence $(\ti y,T_i(e'_j))=(\ti y',T_i(e'_j))=0$, $\ti y+\ti y'=T_i(e'_j)$ where $\ti y=T_i(y),\ti y'=T_i(y')$.
If $T_i(e'_j)=e_h$ for some $h\in[1,2d]$ then $\ti y,\ti y'$ ar joined in $V$, as required.
If this condition is not satisfied then $1<i<D$, $j=i-1$ and $T_i(e'_j)=e_j+e_{j+1}+e_{j+2}$.
We have $(\ti y,e_j+e_{j+1}+e_{j+2})=0$, $\ti y+\ti y'=e_j+e_{j+1}+e_{j+2}$. Since 
$\ti y\in V_i$, we have $(\ti y,e_i)=0$ hence $(\ti y,e_{j+1})=0$ so that
$(\ti y,e_j)=(\ti y,e_{j+2})$. We are in one of the two cases below.

(1) We have $(\ti y,e_j)=(\ti y,e_{j+2})=0$.

(2) We have $(\ti y,e_j)=(\ti y,e_{j+2})=1$.
\nl
In case (1) we consider the four term sequence 
$\ti y,\ti y+e_j,\ti y+e_j+e_{j+2},\ti y+e_j+e_{j+1}+e_{j+2}=\ti y'$; any two
consecutive terms of this sequence are joined in the graph $V$.
In case (2) we consider the four term sequence $\ti y,\ti y+e_{j+1},\ti y+e_j+e_{j+1},
\ti y+e_j+e_{j+1}+e_{j+2}=\ti y'$; any two
consecutive terms of this sequence are joined in the graph $V$. We see that in both cases $\ti y,\ti y'$ 
are in the same connected component of $V$; (a) is proved.

\mpb

Let $V_0=\{x\in V;u(x)=0\}$. Note that $0\in V_0$. We show:

(b) {\it If $x\in V_0$ then $x,0$ are in the same component of the graph $V$.}
\nl
We argue by induction on $D$. If $D=0$ there is nothing to prove. Assume now that $D\ge2$. 
If $(x,e_i)=1$ for all $i\in[1,D]$ then 
$$x=e_{[2,3]}+e_{[6,7]}+e_{[10,11]}+\do+e_{[D-2,D-1]}\text{ if $D/2$ is even},$$
$$x=e_{[1,2]}+e_{[5,6]}+e_{[9,10]}+\do+e_{[D-1,D]}\text{ if $D/2$ is odd}.$$
In both cases we have $u(x)\ne0$ contradicting our assumption. Thus we have $(x,e_i)=0$ for some
$i\in[1,D]$. By 1.9(a) we have $x=T_i(x')+ce_i$ for some $x'\in V'$ and some $c\in\FF_2$. 
By 1.11(b) we have $u'(x')=0$.
By the induction hypothesis $x',0$ are in the same connected component of $V'$.
By (a), $T_i(x'),0$ are in the same connected component of $V$. 
Clearly $x,T_i(x')$ are joined in the graph $V$. Hence $x,0$ are joined in the graph $V$. 
We see that (b) holds.

We show:

(c) {\it $V_0$ is a connected component of the graph $V$.}
\nl
If $x,x'$ in $V$ are in the same connected component of $V$ then $u(x)=u(x')$. (We can assume that
$x,x'$ are joined in the graph $V$. Then for some $i\in[1,D]$ 
we have $x=T_i(y)+ce_i$, $x'=T_i(y)+c'e_i$ where $y\in V'$, $c\in\FF_2,c'\in\FF_2$.
By 1.11(b) we have $u(x)=u'(y)$, $u(x')=u'(y)$, hence  $u(x)=u(x')$, as required).
Thus $V_0$ is a union of connected components of $V$. On the other hand, by (b), $V_0$ is contained in a 
connected component of the graph $V$. This proves (c).

\subhead 1.13\endsubhead
We show:

(a) {\it If $B\in S_D$ then $\la B\ra\sub V_0$.}
\nl
We argue by induction on $D$. If $D=0$ there is nothing to prove. Assume that $D\ge2$. If $B=\emp$ there is
nothing to prove. Assume that $B\ne\emp$. We can find $i\in[1,D]$ and $B'\in S_{D-2}$ such that $B=t_i(B')$.
By 1.9(c) we have $\la B\ra=T_i(\la B'\ra)\op\FF_2e_i$. Using 1.11(b), to prove that $u=0$ on $\la B\ra$ it is 
enough to prove that $u'=0$ on $\la B'\ra$ and this follows from the induction hypothesis. This proves (a).

We show:

(b) {\it If $x\in V_0$ then $x\in\la B\ra$ for some $B\in S_d$.}
\nl
We argue by induction on $D$.  If $D=0$ there is nothing to prove. Assume that $D\ge2$. 
As in the proof of 1.12(b), from the fact that $u(x)=0$ we can deduce that $(x,e_i)=0$ for some $i\in[1,D]$.
By 1.9(a) we have $x=T_i(x')+ce_i$ for some $x'\in V'$ and some $c\in\FF_2$. 
By 1.11(b) we have $u'(x')=0$. By the induction hypothesis we have $x'\in\la B'\ra$ for some $B'\in S_{D-2}$.
Then $x\in T_i(\la B'\ra)\op\FF_2e_1=\la B\ra$ (we use 1.9(c)). This proves (b).

\mpb

From (a),(b) we deduce:

(c) {\it We have $\cup_{B\in S_D}\la B\ra=V_0$.}
\nl
A closely related result is proved in \cite{\SYMPL, 3.4}.

\subhead 1.14\endsubhead
The function $\e:S_D@>>>V$ has values in $\cup_{B\in S_D}\la B\ra$ (see 1.9(c)) hence by 1.13(c) it has values
in $V_0$. Thus, it can be viewed as a function $\e:S_D@>>>V_0$.
\nl
From 1.10(a) we see that:

(a) {\it $\e:S_D@>>>V_0$ is injective.}

\subhead 1.15\endsubhead
Let $F_0$ be the $\QQ$-vector space consisting of functions $V_0@>>>\QQ$. For $x\in V_0$ let $\ps_x\in F_0$ be the
characteristic function of $x$. For $B\in S_D$ let $\Ps_B\in F_0$ be the characteristic function of 
$\la B\ra$. (We use that $\la B\ra\sub V_0$, see 1.13.) Let $\tF_0$ be the $\QQ$-subspace of $F_0$ generated by 
$\{\Ps_B;B\in S_D\}$. When $D\ge2$ we define $\ps'_{x'}$ for $x'\in V'$ and $\Ps'_{B'}$ for $B'\in S_{D-2}$, 
$F'_0,\tF'_0$, in terms of $S_{D-2}$ in the same way as $\ps_x,\Ps_B,F_0,\tF_0$ were defined in terms of $S_D$. 
For any $i\in[1,D]$ we define a linear map $\th_i:F'_0@>>>F_0$ by $f'\m f$ where 
$f(T_i(x')+ce_i)=f'(x')$ for $x'\in V',c\in\FF_2$, $f(x)=0$ for $x\in V-e_i^\pe$. We have 

$\th_i(\ps'_{x'})=\ps_{T_i(x')}+\ps_{T_i(x')+e_i}$ for any $x'\in V'$,

$\th_i(\Ps'_{B'})=\Ps_{t_i(B')}$ for any $B'\in S_{D-2}$. 
\nl
We show:

(a) {\it For any $x\in V_0$, we have $\ps_x\in\tF_0$.}
\nl
We argue by induction on $D$. If $D=0$ the result is obvious. We now assume that $D\ge2$. We first show:

(b) {\it If $x,\tx$ in $V_0$ are joined in the graph $V$ and if (a) holds for $x$, then (a) holds for $\tx$.}
\nl
We can find $j\in[1,2d]$ such that $x+\tx=e_j$, $(x,e_j)=0$.
We have $x=T_j(x')+ce_j$, $\tx=T_j(x')+c'e_j$ where $x'\in V'$ and $c\in\FF_2,c'\in\FF_2$, $c+c'=1$.
 By the induction 
hypothesis we have $\ps'_{x'}=\sum_{B'\in S_{D-2}}a_{B'}\Ps'_{B'}$ where $a_{B'}\in\QQ$. 
Applying $\th_j$ we obtain
$$\ps_x+\ps_{\tx}=\sum_{B'\in S_{D-2}}a_{B'}\Ps_{t_j(B')}$$ 
We see that $\ps_x+\ps_{\tx}\in\tF_0$. Since $\ps_x\in\tF$, by assumption, we see that $\ps_{\tx}\in\tF$. 
This proves (b).

We now prove (a).
Since $V_0$ is the connected component of $V$ containing $0$, 
to prove (a) it is enough (by (b)) to show that (a) holds when $x=0$. This follows from the fact that 
$\ps_0=\Ps_B$ where $B=\emp$. This proves (a).

Since $\tF_0\sub F_0$, we see that (a) implies:

(c) $F_0=\tF_0$.
\nl
We have the following result.
\proclaim{Theorem 1.16}(a) $\{\Ps_B;B\in S_D\}$ is a $\QQ$-basis of $F_0$.

(b) $\e:S_D@>>>V_0$ is a bijection.
\endproclaim
From the definition of $\tF_0$ we have $\dim\tF_0\le|S_D|$. 
By 1.14(a) we have $|S_D|\le|V_0|=\dim F_0$. Since $F_0=\tF_0$ (see 1.15(c)) it follows that 
$\dim\tF_0=|S_D|=|V_0|=\dim F_0$. Using again the definition of $\tF_0$ and the equality $F_0=\tF_0$ we see 
that (a) holds. Since the map in (b) is injective (see 1.14(a))
and $|S_D|=|V_0|$ we see that it is a bijection so that (b) holds.

\subhead 1.17\endsubhead
In this subsection we describe the bijection in 1.16(b) assuming that $D$ is  $2$, $4$ or $6$.
In each case we give a table in which there is one row for each $B\in S_D$; the row corresponding to $B$
is of the form $<B>:(\do)$ where $B$ is represented by the list of intervals of $B$ (we write an interval
such as $[4,6]$ as $456$) and $(\do)$ is a list of the vectors in $\la B\ra$ (we write
$1235$ instead of $e_1+e_2+e_3+e_5$, etc). In each list $(\do)$ we single out the vector corresponding $\e(B)$
in 1.16(b) by putting it in a box. Any non-boxed entry in $(\do)$ appears as a boxed entry in some previous row.
We see that in these cases, 1.5(i) holds.

\mpb

The table for $D=2$.

\mpb

$\emp:(\bx{0})$

$<1>:(0,\bx{1})$

$<2>:(0,\bx{2})$.

\mpb

The table for $D=4$.

$\emp:(\bx{0})$

$<1>:(0,\bx{1})$

$<2>:(0,\bx{2})$

$<3>:(0,\bx{3})$

$<4>:(0,\bx{4})$

$<1,3>:(0,1,3,\bx{13})$

$<1,4>:(0,1,4,\bx{14})$

$<2,4>:(0,2,4,\bx{24})$ 

$<2,123>:(0,2,13,\bx{123})$

$<3,234>:(0,3,24,\bx{234})$.

\mpb

The table for $D=6$.

\mpb

$\emp:(\bx{0})$

$<1>:(0,\bx{1})$

$<2>:(0,\bx{2})$

$<3>:(0,\bx{3})$

$<4>:(0,\bx{4})$

$<5>:(0,\bx{5})$

$<6>:(0,\bx{6})$

$<1,4>:(0,1,4,\bx{14})$

$<1,6>:(0,1,6,\bx{16})$

$<2,4>:(0,2,4,\bx{24})$

$<2,5>:(0,2,5,\bx{25})$

$<2,6>:(0,2,6,\bx{26})$

$<3,6>:(0,3,6,\bx{36})$

$<4,6>:(0,4,6,\bx{46})$

$<1,3>:(0,1,3,\bx{13})$

$<1,5>:(0,1,5,\bx{15})$

$<3,5>:(0,3,5,\bx{35})$

$<2,123>:(0,2,13,\bx{123})$

$<3,234>:(0,3,24,\bx{234})$

$<4,345>:(0,4,35,\bx{345})$

$<5,456>:(0,5,46,\bx{456})$

$<1,3,5>:(0,1,3,5,13,15,35,\bx{135})$

$<1,3,6>:(0,1,3,6,13,16,36,\bx{136})$

$<1,4,345>:(0,1,4,345,14,35,135,\bx{1345})$

$<1,4,6>:(0,1,4,6,14,16,46,\bx{146})$

$<2,4,6>:(0,2,4,6,24,26,46,\bx{246})$

$<1,5,456>:(0,1,5,456,15,46,146,\bx{1456})$

$<2,5,456>:(0,2,5,456,25,46,246,\bx{2456})$

$<2,5,123>:(0,2,5,123,25,13,135,\bx{1235})$

$<2,6,123>:(0,2,6,123,26,13,136,\bx{1236})$

$<2,4,12345>:(0,2,4,24,1345,1235,135,\bx{12345})$    

$<3,234,12345>: (0,3,234,12345,24,15,135,\bx{1245})$

$<3,6,234>:(0,3,6,234,24,36,246,\bx{2346})$

$<3,5,23456>:(0,3,5,2456,35,2346,246,\bx{23456})$

$<4,345,23456>:(0,4,345,23456,35,26,246,\bx{2356})$.

\head 2. The sets $\cf_*(V),\cf(V)$\endhead
\subhead 2.1\endsubhead
We no longer assume that $D$ is even.
We define a collection $\cf_*(V)$ and a collection $\cf(V)$ of subspaces of $V$ by induction on $D$ as follows. 
If $D\in\{0,1\}$, $\cf_*(V)$ and $\cf(V)$ consist of $\{0\}$. If $D\ge2$, a subspace $X$ of $V$ is said to be in 
$\cf_*(V)$ if there exists $i\in[1,D]$ (if $D$ is even) or $i\in[1,D-1]$ (if $D$ is odd) and $X'\in\cf_*(V')$ 
such that $X=T_i(X')\op\FF_2e_i$; a subspace $X$ of $V$ 
is said to be in $\cf(V)$ if either $X=0$ or if there exists $i\in[1,D]$ (if $D$ is even) or
$i\in[1,D-1]$ (if $D$ is odd) and $X'\in\cf(V')$ such that $X=T_i(X')\op\FF_2e_i$. By induction on $D$
we see that for $X\in\cf_*(V)$ we have $X\in\cf(V)$ and $\dim(X)=D/2$ if $D$ is even, $\dim(X)=(D-1)/2$ if $D$ 
is odd. When $D$ is odd, let $\un V$ be the subspace of $V$ with basis $\{e_1,e_2,\do,e_{D-1}\}$. This 
vector space with basis is of the same kind as $V$ in 1.9 (but of even dimension) hence $\cf(\un V),\cf_*(\un V)$
are defined. Using induction on $D$ we see that for $D$ odd we have $\cf(V)=\cf(\un V),\cf_*(V)=\cf_*(\un V)$.
Thus, the study of $\cf(V),\cf_*(V)$ when $D$ is odd is reduced to the similar study when $D$ is even.

We now assume that $D$ is even.
If $B\in S_D$ then $\la B\ra\in\cf(V)$ (this follows from
1.9(c) by induction on $D$). Conversely, if $X\in\cf(V)$ then there exists $B\in S_D$ such that $X=\la B\ra$ (this
again follows from 1.9(c) by induction on $D$). Thus we have a surjective map $S_D@>>>\cf(V)$, $B\m\la B\ra$. We 
show:

(a) {\it This map is a bijection.}
\nl
Indeed, if $B,\tB$ in $S_D$ satisfy $\la B\ra=\la\tB\ra$ then the functions $\Ps_B,\Ps_{\tB}$ in $F_0$
coincide hence $B=\tB$ by 1.16(a). This proves (a).

For $B\in S_D$ we show:

(b) {\it $\{e_I;I\in B\}$ is an $\FF_2$-basis of $\la B\ra$.}
\nl
We argue by induction on $D$. If $D=0$ there is nothing to prove. Assume that $D\ge2$. If $B=\emp$ then (b) is 
obvious. We now assume that $B\ne\emp$. Assume that $\sum_{I\in B}c_Ie_I=0$ with $c_I\in\FF_2$ not 
all zero.
We can find $I=[a,b]\in B$ with $c_I\ne0$ and $|I|$ maximal. If $I'\in B$ is such that $a\in I'$,
$I'\ne I$, $c_{I'}\ne0$, then by $(P_0)$ we have $I\prec I'$ (contradicting the
maximality of $|I|$) or $I'\prec I$ (contradicting $a\in I'$). Thus no $I'$ as above exists.
Thus when $\sum_{I_1\in B}c_{I_1}e_{I_1}$ is written in the basis $\{e_j;j\in[1,D]\}$, the coefficient of $e_a$
is $c_{I_1}$ hence $c_{I_1}=0$, contradicting $c_{I_1}\ne0$. This proves (b).

We show:

(c) {\it If $X\in\cf(V)$ then $X$ is an isotropic subspace of $V$.}
\nl
We argue by induction on $D$. If $D=0$ there is nothing to prove. Assume that $D\ge2$. If $X=0$ then (c) is 
obvious. We now assume that $X\ne0$. Then there exists $i\in[1,D]$ and $X'\in\cf(V')$ such that 
$X=T_i(X')\op\FF_2e_i$. By the induction hypothesis, $X'$ is isotropic in $V'$. Since $T_i$ is compatible with
the symplectic forms it follows that $T_i(X')$ is an isotropic subspace of $V$. Since $T_i(X')$ is contained in
$e_i^\pe$, $T_i(X')\op\FF_2e_i$ is also isotropic. This proves (c). Alternatively, (c) can be deduced from
property $(P_0)$.

\subhead 2.2\endsubhead
For $\d\in\{0,1\}$ let $[1,D]^\d=\{i\in[1,D];i=_2\d\}$. 
Let $V^\d$ be the subspace of $V$ with basis $\{e_i;i\in[1,D]^\d\}$. We have $V=V^0\op V^1$.
Similarly, if $D\ge2$,  we have $V'=V'{}^0\op V'{}^1$ where $V'{}^\d$ has basis $\{e'_i;i\in[1,D-2]^\d\}$. 

For any $I\in\ci^1_D$ and $\d\in\{0,1\}$ we set $I^\d=I\cap[1,D]^\d$, so that $I=I^0\sqc I^1$; 
we define $\k(I)\in\{0,1\}$ by $a=_2\k(I)$ or equivalently $b=_2\k(I)$ where $I=[a,b]$. We show:

(a) {\it Let $B\in S_D$ and let $I\in B$. Let $\d=\k(I)$. We have $e_{I^\d}=\sum_{I'\in B;I'\sub I}e_{I'}$.}
\nl
We argue by induction on $|I|$. If $|I|=1$ the result is obvious. Assume now that $|I|>1$. By $(P_0),(P_1)$, we 
can find $[a_1,b_1],[a_2,b_2],\do,[a_k,b_k]$ in $B$ such that $a_1\le b_1<a_2\le b_2<a_3\le b_3<\do$, 
$a_1,b_1,a_2,b_2,\do$ are all in $1-\d+2\ZZ$  
and $[a,b]\cap(1-\d+2\ZZ)\sub[a_1,b_1]\cup[a_2,b_2]\cup\do\cup[a_k,b_k]$. From the definition we have
$e_{I^\d}=e_I+\sum_{j=1}^ke_{[a_j,b_j]^{1-\d}}$. By the induction hypothesis, for $j\in[1,k]$ we have
$e_{[a_j,b_j]^{1-\d}})=\sum_{I'\in B;I'\sub[a_j,b_j]}e_{I'}$. Thus we have 
$$e_{I^\d}=e_I+\sum_{I'\in B;I'\sub\cup_j[a_j,b_j]}e_{I'}=\sum_{I'\in B;I'\sub I}e_{I'}.$$
This proves (a).

We show:

(b) {\it Let $B\in S_D$. Then $\{e_{I^{\k(I)}};I\in B\}$ is a basis of the vector space $\la B\ra$.}
\nl
From (a) we see that the collection of vectors $\{e_{I^{\k(I)}};I\in B\}$ is related to the collection of vectors
$\{e_I;I\in B\}$ by an upper triangular matrix with $1$ on diagonal. Hence the result follows from 2.1(b).

We deduce that if $B\in S_D$ and $X=\la B\ra\in\cf(V)$ then for $\d\in\{0,1\}$,

(c) {\it $X^\d=X\cap V^\d$ has basis $\{e_{I^{\k(I)}};I\in B,\k(I)=\d\}$; in particular, $X=X^0\op X^1$.}

\subhead 2.3\endsubhead
Assume that $D\ge2$. Let $i\in[1,D]$ and let $\d\in\{0,1\}$. There is a unique linear map
$T_i^\d:V'{}^\d@>>>V^\d$ such that

$T_i^\d(e'_k)=e_k$ if $k\le i-2$, $k=_2\d$;

$T_i^\d(e'_k)=e_{k+2}$  if $k\ge i$, $k=_2\d$;

$T_i^\d(e'_{i-1})=e_{i-1}+e_{i+1}$ if $i=_2\d+1$, $1<i<D$.

Note that $T_i^\d$ is injective and $(x,y)'=(T_i^0(x),T_i^1(y))$ for any $x\in V'{}^0,y\in V'{}^1$.
For any $I'\in\ci^1_{D-2}$ such that $\k(I')=\d$ we have $T_i^\d(e'_{I'{}^\d})=e_{\x_i(I')^\d}$.
(Here $\k(I'),I'{}^\d$ are defined in terms of $I'$ in the same way as $\k(I),I^\d$ are defined in 2.2.) 
Let $V_i^\d$ be the image of $T_i^\d:V'{}^\d@>>>V^\d$. From the definitions we deduce:

(a) {\it  We have $V_i\op\FF_2e_i=V_i^0\op V_i^1\op\FF_2e_i$.}
\nl
We define a collection $\cc(V^\d)$ of subspaces of $V^\d$ by induction on $D$ as follows. If $D=0$, $\cc(V^\d)$ 
consists of $\{0\}$. If $D\ge2$, a subspace $\cl$ of $V^\d$ is said to be in $\cc(V^\d)$ if either $\cl=0$ or if 
there exists $i\in[1,D]$ and $\cl'\in\cc(V'{}^\d)$ such that $\cl=T_i^\d(\cl')\op\FF_2e_i$ (if 
$i=_2\d$) or 
$\cl=T_i^\d(\cl')$ (if $i=_2\d+1$).

We show:

(b) {\it If $X\in\cf(V)$ then $X^\d\in\cc(V^\d)$.}
\nl
We argue by induction on $D$. If $D=0$ the result is obvious. Assume now that $D\ge2$. If $X=0$ there  is
nothing to prove. Assume that $X\ne0$. We can find $i\in[1,D]$ and $X'\in\cf(V')$ such that 
$X=T_i(X')\op\FF_2e_i$. By the induction hypothesis we have $X'{}^\d\in\cc(V'{}^\d)$. Hence
$T_i^\d(X'{}^\d)\op\FF_2e_i\in\cc(V^\d)$ if $i=_2\d$, $T_i^\d(X'{}^\d)\in\cc(V^\d)$ if $i=_2\d+1$. It is 
enough to prove that $T_i^\d(X'{}^\d)\op\FF_2e_i=X^\d$ if $i=_2\d$, 
$T_i^\d(X'{}^\d)=X^\d$ if $i=_2\d+1$,
or that $T_i^\d(X'{}^\d)\op\FF_2e_i=(T_i(X')\op\FF_2e_i)\cap V^\d$ if $i=_2\d$, 
$T_i^\d(X'{}^\d)=(T_i(X')\op\FF_2e_i)\cap V^\d$ if $i=_2\d+1$.
This follows by comparing the definition of $T_i^\d$ with that of $T_i$.

\subhead 2.4\endsubhead
Let $\d\in\{0,1\}$. If $Z$ is a subspace of $V^\d$ we set $Z^!=\{x\in V^{1-\d}; (x,z)=0\qua\frl z\in Z\}$. 
Similarly, if $Z'$ is a subspace of $V'{}^\d$ we set $Z'{}^!=\{x\in V'{}^{1-\d}; (x,z)'=0\qua\frl z\in Z'\}$. 
Let $\cl\in\cc(V^\d)$. We show:

(a) {\it We have $\cl^!\in\cc(V^{1-\d})$ and $\cl\op\cl^!\sub V$ is in $\cf(V)$.}
\nl
The first statement of (a) follows from the second statement, using 2.3(b). We prove the second statement of (a) 
by induction on $D$. If $D=0$ the result is immediate. Assume now that $D\ge2$. If $\cl=0$, then 
$\cl^!=V^{1-\d}=\la B\ra$ where $B=\{\{j\};j\in[1,D]^{1-\d}\}\in S_D$; thus we have $\cl^!\in\cf(V)$. Next we assume
that $\cl\ne0$. We can find $i\in[1,D]$ and $\cl'\in\cc(V'{}^\d)$ such that 
$\cl=T_i^\d(\cl')\op\FF_2e_i$ 
(if $i=_2\d$) or $\cl=T_i^\d(\cl')$ (if $i=_2\d+1$). 
By the induction hypothesis we have $\cl'\op\cl'{}^!\in\cf(V')$. Hence $T_i(\cl'\op\cl'{}^!)\op\FF_2e_i\in\cf(V)$.
From the definition we have $T_i(\cl'\op\cl'{}^!)\op\FF_2e_i
=T_i^\d(\cl')\op T_i^{1-\d}(\cl'{}^!)\op\FF_2e_i$. 
Thus we have $T_i^\d(\cl')\op T_i^{1-\d}(\cl'{}^!)\op\FF_2e_i\in\cf(V)$ or equivalently 
$\cl\op T_i^{1-\d}(\cl'{}^!)\in\cf(V)$ (if $i=_2\d$) and $\cl\op T_i^{1-\d}(\cl'{}^!)\op\FF_2e_i\in\cf(V)$ 
(if $i=_2\d+1$). It is enough to show: $\cl^!=T_i^{1-\d}(\cl'{}^!)$ if $i=_2\d$ and 
$\cl^!=T_i^{1-\d}(\cl'{}^!)\op\FF_2e_i$ if $i=_2\d+1$. If $y\in\cl'{}^!$, $x\in\cl'$, we have 
$(T_i^{1-\d}(y),T_i^\d(x))=(y,x)'=0$; if $i=_2\d$ we have $(T_i^{1-\d}(y),e_i)=0$. If $i=_2\d+1$ 
we have 
$(e_i,T_i^\d(x))=0$. We see that $T_i^{1-\d}(\cl'{}^!)\sub\cl^!$ 
if $i=_2\d$ and $T_i^{1-\d}(\cl'{}^!)\op\FF_2e_i\sub\cl^!$ if $i=_2\d+1$. The last two inclusions are 
between vector spaces of the same dimension; hence they must be equalities. This completes the proof of (a).

\mpb

Let $S_{D,*}=\{B\in S_D;|B|=D/2\}$. From 2.1(b) we see that the bijection $S_D@>\si>>\cf(V)$, $B\m\la B\ra$ (see 
2.1(a)) restricts to a bijection

(b) $S_{D,*}@>\si>>\cf_*(V)$.
\nl
We show:

(c) {\it We have a bijection $\io:\cc(V^\d)@>\si>>\cf_*(V)$ given by $\io(\cl)=\cl\op\cl^!$.}
\nl
The fact that $\io$ is well defined follows from (a). (For $\cl\in\cc(V^\d)$ we have $\dim(\cl\op\cl^!)=D/2$.)
We define $\io':\cf_*(V)@>>>\cc(V^\d)$ by $X\m X^\d$. This is well defined by 2.3(b). Clearly, $\io'\io=1$. Let 
$X\in\cf_*(V)$. Then $X^{1-\d}\sub(X^\d)^!$ since $X$ is isotropic so that $X^\d\op(X^\d)^!\sub X$; this is an 
inclusion of vector spaces of the same dimension, hence is an equality. Thus $\io\io'=1$. This proves that $\io$ 
is a bijection.

\subhead 2.5\endsubhead
Let $\d\in\{0,1\}$. 
We define a subset $S_D^\d$ of $R^1_D$ by induction on $D$ as follows.
When $D=0$, $S_D^\d$ consists of $\emp\in R^1_D$. When $D\ge2$ we say that $\b\in R^1_D$ is in $S^\d_D$ if either 
$\b=\emp$ or if 

(i) there exists $i\in[1,D]$ and $\b'\in S^\d_{D-2}$ such that $\b=\{\x_i(I');I'\in\b'\}\sqc\{i\}$ if 
$i=_2\d$ and $\b=\{\x_i(I');I'\in\b'\}$ if $i=_2\d+1$.
\nl
From the definition we see by induction on $D$ that if $\b\in S^\d_D$ and $I\in\b$ then 
$\k(I)=\d$.

Let $S'_D{}^\d$ be the set of all $\b\in R^1_D$ such that $\k(I)=\d$ for any $I\in\b$ and such that the
following holds.

$(P^\d_0)$ {\it If $I\in\b$, $\tI\in\b$, then either $I=\tI$, or $I\spa\tI$, or $I\prec\tI$, or $\tI\prec I$.}
\nl
By arguments similar to those in 1.3 we see that 

(a) {\it We have $S^\d_D=S'_D{}^\d$.}
\nl
We show:

(b) {\it If $B\in S_D$ then ${}^\d B:=\{I\in B;\k(I)=\d\}$ is in $S^\d_D$.}
\nl
From 2.5(c) we see that ${}^\d B\in S'_D{}^\d$ hence (using (a)) ${}^\d B\in S_D^\d$. 

\mpb

Using the definitions we can verify:

(c) {\it Assume that $D\ge2$, that $B'\in S_{D-2}$ and that $B=t_i(B')\in S_D$. Let $\b'={}^\d B'\in S^\d_{D-2}$,
$\b={}^\d B\in S^d_D$. Then $\b$ is obtained from $\b'$ as in (i) above.}
\nl
Let ${}'S_D^\d$ be the set of all subsets of $R^1_D$ of the form ${}^\d B$ for some $B\in S_{D,*}$. We show:

(d) {\it ${}'S^\d_D=S^\d_D$.}
\nl
The inclusion ${}'S^\d_D\sub S^\d_D$ follows from (b). Conversely we show that if $\b\in S^\d_D$
then $\b\in{}'S^\d_D$. We argue by induction on $D$. When $D=0$ there is nothing to prove. Assume that $D\ge2$.
If $\b=\emp$ there is nothing to prove. Assume that $\b\ne\emp$. We can find $i\in[1,D]$ and $\b'\in S^\d_{D-2}$
such that $\b$ is obtained from $\b'$ as in (i) above. By the induction hypothesis we have
$\b'={}^\d B'$ where $B'\in S_{D-2,*}$. Let $B=t_i(B')$. We have $B\in S_{D,*}$.
Let $\ti\b={}^\d B\in{}'S^\d_D$. By (c), $\ti\b$ is obtained from $\b'$ as in (i) above. Since $\b$ has the
same property, we have $\ti\b=\b$. Thus $\b\in{}'S^d_D$, as required. This proves (d).

We show:

(e) {\it The map $S_{D,*}@>>>{}'S^\d_D$, $B\m{}^\d B$ is a bijection.}
\nl
It is enough to show that this map is injective. Assume that $B\in S_{D,*},\tB\in S_{D,*}$ satisfy
${}^\d B={}^\d\tB$. We must show that $B=\tB$.
By the proof of 2.4(c) we have a bijection $\io':\cf_*(V)@>>>\cc(V^\d)$ given by $X\m X^\d$.
Now $\io'(\la B\ra)$ has basis $\{e_{I^{\k(I)}};I\in B,\k(I)=\d\}$ and
$\io'(\la\tB\ra)$ has basis $\{e_{I^{\k(I)}};I\in\tB,\k(I)=\d\}$. Since ${}^\d B={}^\d\tB$, these two bases
coincide hence $\io'(\la B\ra)=\io'(\la\tB\ra)$. Since $\io'$ is a bijection we deduce that
$\la B\ra=\la\tB\ra$. Using 2.1(a) we see that $B=\tB$. This proves (e).

\mpb

Combining (d),(e) we obtain:

(f) {\it The map $S_{D*}@>>>S^\d_D$, $B\m{}^\d B$ is a bijection.}
\nl
For any $\b\in S^\d_D$ let $\la\b\ra$ be the $\FF_2$-subspace of $V^\d$ spanned by $\{e_{I^{\k(i)}};I\in\b\}$. 
By the proof of (e), we have $\la\b\ra\in\cc(V^\d)$ and $\dim\la\b\ra=|\b|$. We show:

(g) {\it The map $\b\m\la\b\ra$ is a bijection $\io'':S^\d_D@>\si>>\cc(V^\d)$.}
\nl
We have a commutative diagram
$$\CD
S_{D,*}@>>>            \cf_*(V)\\
@VVV                 @V\io'VV  \\
S^\d_D@>\io''>>\cc(V^\d)    \\
\endCD$$
where the top horizontal map is a bijection as in 2.4(b), the left vertical map is a bijection as in (e) 
(see also (d)) and $\io'$ is a bijection as in the proof of (e). It follows that $\io''$ is a bijection.
This proves (g).

\subhead 2.6\endsubhead
Let $\d\in\{0,1\}$. We define a bijection $S_D^\d@>\si>>S_D^{1-\d}, \b@>>>\b^!$ as follows.
Let $\b\in S_D^\d$. By 2.5(g), 
we have $\la\b\ra\in\cc(V^\d)$ and by 2.4(a) we have $\la\b\ra^!\in\cc(V^{1-\d})$.
Then $\b^!$ is the unique element of $S_D^{1-\d}$ such that $\la\b\ra^!=\la\b^!\ra$, see 2.5(g).
From the definition we have $(\b^!)^!=\b$ and $|\b^!|=(D/2)-|\b|$. 
Recall that $\la\b\ra\op\la\b^!\ra=\la B\ra$ where $B\in S_{D,*}$ satisfies ${}^\d B=\b,{}^{1-\d}B=\b^!$. 

The order reversing involution $i\m i^*=D+1-i$ of $[1,D]$ induces an involution $R^1_D@>>>R^1_D$, 
$I\m I^*=\{i^*;i\in I\}$ and an involution $S_D@>>>S_D$, $B\m B^*:=\{I^*;I\in B\}$.
It also induces a bijection $\g_\d:S^{1-\d}_D@>\si>>S^\d_D$. 
Then $\b\m\g_\d(\b^!)$ is a bijection $S^\d_D@>>>S^\d_D$ which carries any subset with $m$
elements ($m\in[0,D/2]$) to a subset with $(D/2)-m$ elements.

\subhead 2.7\endsubhead
Let $\d\in\{0,1\}$. Let $U^\d=\{(\cl,\cl')\in\cc(V^\d)\T\cc(V^\d);\cl\sub\cl'\}$. We define a map 

(a) $\cf(V)@>>>U^\d$ by $X\m(X^\d,(X^{1-\d})^!)$.
\nl
(We have $X^\d\sub(X^{1-\d})^!$ since $X$ is isotropic.) This map is injective since $X$ can be reconstructed 
from $X^\d,X^{1-\d}$: we have $X=X^\d\op X^{1-\d}$.

We note that the map (a) is not surjective. For example, if $D=2,\d=0$ and $\cl=0$, $\cl'=\FF_2e_2$
then $(\cl,\cl')\in U^0$ is not in the image of the map (a). The following result is a reformulation of 2.4(c).

(b) {\it The map (a) restricts to a bijection $\cf_*(V)@>\si>>\{(\cl,\cl')\in U^\d;\cl=\cl'\}$.}

\subhead 2.8\endsubhead
In the remainder of this section we prove Theorem 0.4 assuming that $W$ is a Weyl group of type $B_n,C_n$ or 
$D_n$.
If $|c|=1$ the theorem is trivial; we have $\cg_c=\{1\}$ and $\BB_c$ consists of the unique 
representation in $c$.  Assume now that $|c|\ge2$. As in \cite{\ORA, 4.5, 4.6}, \cite{\CLASSII}, \cite{\CELL}, 
we can find $D\in\{2,4,6,\do\}$ and $\d\in\{0,1\}$ such that if $V$ is the $\FF_2$-vector 
space with basis $\{e_i;i\in[1,D]\}$ as in 1.9, then (i)-(iii) below hold.

(i) The group $\cg_c$ in 0.3 is $V^\d$; hence $M(\cg_c)=V^\d\op\Hom(V^\d,\CC^*)$ can be identified with 
$V=V^\d\op V^{1-\d}$ (an element $y\in V^{1-\d}$ can be identified with the homomorphism $V^\d@>>>\CC^*$
given by $x\m (-1)^{(x,y)}$).

(ii) $c$ is naturally in bijection with $V_0$ (see 1.12); hence any object $\ce\in\car_c$ can be viewed as
the function $f_\ce:V_0@>>>\NN$ such that for $E\in c$ the multiplicity of $E$ in $\ce$ is equal to the value of
$f_\ce$ at the point of $V_0$ corresponding to $E$;

(iii) The constructible representations in $\car_c$ viewed as functions $V_0@>>>\NN$ are exactly the
characteristic functions of the subsets $X\sub V$ with $X\in\cf_*(V)$.
\nl
(More accurately, the results in {\it loc.cit.} for $W$ of type $D_n$ are formulated in terms
of a $V$ as in 1.9 with odd $D$, but they can be restated in terms of a $V$ as in 1.9 with $D$ even, by
the argument in the first part of 2.1.)

If $\cl$ is a subspace of $V^\d$ then $S_{\cl,\cl}\in\CC[M(\cg_c)]=\CC[V]$ (see (i) and 0.2) can be identified
with the function $V@>>>\CC$ whose value is $1$ at any element of $\cl\op\cl^!$ and is $0$ at any
element of $V-(\cl\op\cl^!)$. If $\cl\in\cc(V^\d)$ this is the characteristic function of some $X\in\cf_*(V)$
namely, $X=\cl\op\cl^!$; the converse also holds. We see that $\fF_c$ (see 0.3) consists 
of the subspaces 
$\cl\in\cc(V^\d)$. We have $0\in\cc(V^\d)$ hence $\hat\fF_c=\fF_c$.
Now $\ti\Th_c$ becomes the set of pairs $(\cl,\cl')\in\cc(V^\d)\T\cc(V^\d)$ such that $\cl\sub\cl'$.
We define $\Th_c$ to be the set of pairs $(\cl,\cl')\in\cc(V^\d)\T\cc(V^\d)$ such that 
$\cl\op\cl'{}^!\in\cf(V)$. (We then automatically have $\cl\sub\cl'$ since the subspaces in $\cf(V)$ are
isotropic. Thus $\Th_c\sub\ti\Th_c$.) If $(\cl,\cl')\in\ti\Th_c$ then 
$S_{\cl,\cl'}\in\CC[M(\cg_c)]=\CC[V]$ (see (i) and 0.2) can be identified
with the function $V@>>>\CC$ whose value is $1$ at any element of $\cl\op\cl'{}^!$ and is $0$ at any
element of $V-(\cl\op\cl'{}^!)$. If $(\cl,\cl')\in\Th_c$, this is the characteristic function of some 
$X\in\cf(V)$ namely $X=\cl\op\cl'{}^!$; the converse also holds. We see that $\Th_c$ can be identified
with $\cf(V)$. With these identification Theorem 0.4 follows from the results in \S1,\S2. The representations
in $\BB_c$ corespond as in (ii) to the functions $f^X:V_0@>>>\NN$ which equal $1$ on $X$ and equal $0$ on $V_0-X$
(where $X\in\cf(V)$). The bijection $c@>>>\BB_c$ mentioned in 0.5 is $x\m\la\e\i(x)\ra$ where $\e$ is as in 
1.16(b).

\head 3. Exceptional Weyl groups\endhead
\subhead 3.1\endsubhead
In this section we will prove Theorem 0.4 assuming that $W$ is of exceptional type.
In 3.2-3.8 we will give a table of new representations in $\car_c$ in the form a matrix $M_c$
indexed by $c\T c$. (The table will be justified in 3.10.)
The columns of $M_c$ are indexed by the representations in $c$. 
The rows of $M_c$ are also indexed by the  representations in $c$ (for any $k\in[1,|c|]$, the $k$-th row from up 
to down is indexed by the same representation in $c$ as the $k$-th column from left to right). Each row of $M_c$ 
corresponds to a new representation; the entries of that row give the multiplicities of the various representations in $c$ 
in the new representation. The first row in $M_c$ stands for the special representation in $c$.

\subhead 3.2\endsubhead
If $|c|=1$, $M_c$ is the $1\T1$ matrix with entry $1$.

\subhead 3.3\endsubhead
If $|c|=2$ (so that $W$ is of type $E_7$ or $E_8$) we order $c$ using its bijection with
$\{(1,1),(1,\e)\}$ in \cite{\ORA, 4.12, 4.13} (ordered from left to right); then $M_c$ is
$$\left(\matrix
1 &0  \\
1 & 1  \\
\endmatrix\right)$$
The second row stands for a  constructible representation.

\subhead 3.4\endsubhead
If $|c|=3$ we order $c$ using its bijection with
$\{(1,1),(g_2,1),(1,\e)\}$ in \cite{\ORA, 4.10, 4.11, 4.12, 4.13} (ordered from left to right);
 then $M_c$ is
$$\left(\matrix
1 &0 &0  \\
1 & 1 &0\\
1 & 0 &1\\
\endmatrix\right)$$
The last two rows stand for  constructible representations.

\subhead 3.5\endsubhead
If $|c|=4$ (so that $W$ is of type $G_2$) we order $c$ using its bijection with
$\{(1,1),(1,r),(g_2,1),(g_3,1)\}$ in \cite{\ORA, 4.8} (ordered from left to right); then $M_c$ is
$$\left(\matrix
1 & 0  &0     &0    \\
1 & 1  &0     &0    \\
 1 & 1  &1     &0    \\
1 & 0  &1     &1    \\
\endmatrix\right)$$
The last two rows stand for constructible representations.

\subhead 3.6\endsubhead
If $|c|=5$ (so that $W$ is of type $E_6,E_7$ or $E_8$) we order $c$ using its bijection with
$\{(1,1),(1,r),(g_2,1),(g_3,1),(1,\e)\}$ in \cite{\ORA, 4.11, 4.12, 4.13} (ordered from left to right);
then $M_c$ is
$$\left(\matrix
1 &0      & 0   & 0   &0 \\
  1& 1   & 0   & 0   &0 \\
 1 &   1     &1  & 0   &0 \\
 1 &  0     & 1   &1 &0 \\
 1 &  2     & 0   & 0      &1 \\
\endmatrix\right)$$
The last three rows stand for constructible representations.

\subhead 3.7\endsubhead
If $|c|=11$ (so that $W$ is of type $F_4$) we write the elements of $c$
(notation of \cite{\ORA, 4.10}) in the order
$$12_1,9_3,6_2,1_3,16_1,9_2,4_4,6_1,4_3,4_1,1_2$$
(from left to right); then $M_c$ is:
$$\left(\matrix
1  &0    &0   &0    &0   &0      &0   &0    &0  &0  &0                   \\
  1   &1  &0   &0    &0    &0     &0   &0    &0  &0   &0                 \\
 1   &1    &1  &0    &0    &0     &0   &0    &0  &0   &0               \\
 1    & 2  & 1   &1  & 0    &0    & 0   &0   & 0  & 0   &0                    \\
 1    & 1  & 1   &0   &1     &0   &0   &0    &0  &0   &0             \\
  1   & 0   &1   &0   & 1   &1    &0   &0    &0  &0   &0               \\
 1    & 2   &1   &1    &1    &0   &1   &0    &0  &0   &0               \\
 1    & 1   &0   & 0    &1    &0   & 1   &1   &0  &0   &0                \\
 1    & 0   &0   & 0    &1    &1   & 0    &1   &1 &0     &0            \\
 1    & 1   &1   & 0    &2    &1    &0    &0   & 0  &1   &0              \\
 1    & 0   &1    &0    &1    &2    &0    &0   & 1  &0   &1           \\
\endmatrix\right)$$
The last five rows stand for constructible representations.

\subhead 3.8\endsubhead
If $|c|=17$ (so that $W$ is of type $E_8$) we write the elements of $c$ 
(with notation of \cite{\ORA, 4.13.2} with subscripts omitted) in the order
$$\align&4480,5670,4536,1680,1400,70,7168,5600,3150,4200,2688,2016,\\&448,1134,1344,420,168\endalign$$
(from left to right); then $M_c$ is:
$$\left(\matrix
1  &0   &0    &0   &0   &0   &0   &0      &0  &0    &0   &0   &0   &0    &0   &0                 &0        \\
 1    &1 &0    &0   &0   &0   &0  &0     &0  &0     &0  &0  &0   &0    &0   &0            &0                  \\
 1   &1    &1  &0   &0   &0   &0   &0      &0  &0    &0  &0  &0   &0   &0   &0              &0                \\
 1    &2   &1   &1  &0   &0   &0   &0      &0  &0     &0 &0   &0   &0   &0   &0              &0               \\
 1    &2   &2    &1  &1  &0   &0   &0      &0  &0     &0 &0   &0   &0    &0   &0     &0               \\
 1    &3   &3    &3   &2  &1  &0   &0      &0  &0     &0  &0  &0   &0    &0   &0     &0               \\
 1    &1   &1    &0   &0   &0  &1  &0      &0  &0     &0  &0   &0   &0   &0   &0       &0            \\
 1    &2   &2    &1   &1   &0   &1  &1     &0  &0     &0  &0   &0   &0   &0   &0              &0          \\
 1    &1   &1    &0   &0   &0   &1   &1     &1&0      &0  &0   &0   &0    &0   &0               &0        \\
 1    &1   &1    &0   &1   &0   &1   &1      &0  &1   &0  &0   &0   &0    &0   &0              &0            \\
 1    &2   &2    &1   &1   &0   &2   &2      &0  &1   &1  &0   &0   &0    &0   &0               &0           \\
 1    &1   &1    &0   &0   &0   &2   &1      &1  &1    &1   &1 &0   &0    &0   &0                 &0      \\
 1    &3   &3    &3   &2   &1   &1   &2      &0  &0    &0   &0  &1  &0    &0   &0                &0       \\
 1    &2   &1    &1   &0   &0   &1   &2      &1  &0    &0   &0  &1   &1  &0   &0                &0         \\
 1    &1   &0    &0   &0   &0   &1   &1      &1  &1    &0   &0   &0   &1  &1    &0                 &0     \\
 1    &0   &0    &0   &0   &0   &1   &0      &1  &1    &0   &1   &0   &0  &1   &1        &0        \\
 1    &1   &1    &0   &1   &0   &1   &1      &0  &2    &0    &0   &0   &0  &1   &0     &1  \\
\endmatrix\right)$$
The last seven rows stand for constructible representations.

\subhead 3.9\endsubhead
For $N\ge1$ let $S_N$ be the group of all permutations of $[1,N]$.
If $a_1\ge a_2\ge\do$ is a partition of $N$ (written as $a_1a_2\do$) we say that a subgroup $H$ of $S_N$ is in
$\cs_{a_1a_2\do}$ if $H$ is conjugate to the subgroup of all permutations of $[1,N]$ which keep stable each of
the subsets $[1,a_1], [a_1+1,a_1+a_2],[a_1+a_2+1,a_1+a_2+a_3],\do$.
We say that a subgroup $H$ of $S_N$ (with $N\ge4$) is in $\ti\cs_N$ if it is conjugate to the subgroup of
of all permutations of $[1,N]$ which act as identity on $[1,N]-[1,4]$ and whose restriction to
$[1,4]$ commute with the permutation $1\m4\m1,2\m3\m2$.

The following results come from \cite{\LEA}.

If $|c|=1$ we have $\cg_c=\{1\}$ and $\hat\fF_c$ consists of $\{1\}$. 

In the setup of 3.3 or 3.4 we have $\cg_c=S_2$ and $\hat\fF_c$ consists of $S_2,\{1\}$.

In the setup of 3.5 or 3.6 we have $\cg_c=S_3$ and $\hat\fF_c$ consists of $S_3,\{1\}$ and the subgroups of 
$S_3$ in $\cs_{21}$.

In the setup of 3.7 we have $\cg_c=S_4$ and $\hat\fF_c$ consists of $S_4,\{1\}$ and the
subgroups of $S_4$ in $\cs_{31},\cs_{22},\cs_{211},\ti\cs_4$.

In the setup of 3.8 we have $\cg_c=S_5$ and $\hat\fF_c$ consists of $S_5,\{1\}$ and the
subgroups of $S_5$ in $\cs_{41},\cs_{32},\cs_{311},\cs_{221},\cs_{2111},\ti\cs_5$.

\subhead 3.10\endsubhead
We describe the set $\ti\Th_c$ in each of the cases 3.2-3.8.

If $|c|=1$, $\ti\Th_c$ consists of $(1,1)$. (We shall write $1$ instead of $\{1\}$.)

In the setup of 3.3 or 3.4, $\ti\Th_c$ consists of $(1,S_2),(S_2,S_2),(1,1)$.

In the setup of 3.5 or 3.6, $\ti\Th_c$ consists of $(1,S_3),(1,H_{21}),(H_{21},H_{21}),(S_3,S_3),(1,1)$
where $H_{21}$ runs through $\cs_{21}$.

In the setup of 3.7, $\ti\Th_c$ consists of 
$$\align&(1,S_4),(1,H_{31}),(1,H_{22}),(1,H_{211}),(\ti H_{211},H_{22}),(\ti H_{22},\ti H),\\&
(H_{211},H_{211}),(H_{31},H_{31}),(S_4,S_4),(H_{22},H_{22}),(\ti H,\ti H), (1,1),(1,\ti H)\endalign$$
where $H_{211}$ runs through $\cs_{211}$, $H_{31}$ runs through $\cs_{31}$, $H_{22}$ runs through $\cs_{22}$, 
$\ti H$ runs through $\ti\cs_4$; for $H_{22}\in\cs_{22}$, $\ti H_{211}$ denotes
one of the two subgroups in $\cs_{211}$ contained in $H_{22}$; for $\ti H\in\ti\cs_4$, 
$\ti H_{22}$ denotes the
unique subgroup in $\cs_{22}$ contained in $\ti H$.

In the setup of 3.8, $\ti\Th_c$ consists of 
$$\align&(1,S_5),(1,H_{41}),(1,H_{32}),(1,H_{311}),(1,H_{221}),(1,H_{2111}),(\ti H_{2111},H_{32}),\\&
(\ti H_{2111},H_{221}), (\ti H_{311},H_{32}), (\ti H_{221},\ti H),(H_{221},H_{221}),(H_{32},H_{32}),\\&
(H_{2111},H_{2111}), (H_{311},H_{311}),(H_{41},H_{41}),(S_5,S_5),(\ti H,\ti H),(1,1),(1,\ti H),
\endalign$$
where $H_{2111}$ runs through $\cs_{2111}$, $H_{221}$ runs through $\cs_{221}$, 
$H_{32}$ runs through $\cs_{32}$, $H_{311}$ runs through $\cs_{311}$, $H_{41}$ runs through $\cs_{41}$, 
$\ti H$ runs through $\ti\cs_5$; for $H_{221}\in\cs_{221}$, $\ti H_{2111}$ denotes one of the two subgroups
in $\cs_{2111}$ contained in $H_{221}$; for $H_{32}\in\cs_{32}$, $\ti H_{2111}$ denotes the unique subgroup
in $\cs_{2111}$ which is a normal subgroup of $H_{32}$ and $\ti H_{311}$ denotes the unique subgroup
in $\cs_{311}$ which is a normal subgroup of $H_{32}$; for $\ti H\in\ti\cs_5$, $\ti H_{221}$
denotes the unique subgroup in $\cs_{221}$ contained in $\ti H$.

\subhead 3.11\endsubhead
We define the set $\Th_c$ in each of the cases 3.2-3.8 by removing from $\ti\Th_c$ the pair $(1,1)$
whenever $c$ is anomalous (see 0.3)
and by removing the pairs $(1,\ti H)$ with $\ti H$ in $\ti\cs_4$ or $\ti\cs_5$ whenever $\ti\cs_4$ 
or $\ti\cs_5$ is defined. 
This guarantees that for $(H,H')\in\Th_c$, $H'/H$ is isomorphic to a product of symmetric groups.

If $|c|=1$, $\Th_c$ consists of $(1,1)$.

In the setup of 3.3, $\Th_c$ consists of $(1,S_2),(S_2,S_2)$.

In the setup of 3.4, $\Th_c=\ti\Th_c$ consists of $(1,S_2),(S_2,S_2),(1,1)$.

In the setup of 3.5, $\Th_c$ consists of $(1,S_3),(1,H_{21}),(H_{21},H_{21}),(S_3,S_3)$ (notation of 
3.10).

In the setup of 3.6, $\Th_c=\ti\Th_c$ consists of $(1,S_3),(1,H_{21}),(H_{21},H_{21}),(S_3,S_3),(1,1)$ 
(notation of 3.10).

In the setup of 3.7, $\Th_c$ consists of 
$$\align&(1,S_4),(1,H_{31}),(1,H_{22}),(1,H_{211}),(\ti H_{211},H_{22}),(\ti H_{22},\ti H),\\&
(H_{211},H_{211}),(H_{31},H_{31}),(S_4,S_4),(H_{22},H_{22}),(\ti H,\ti H),\endalign$$
(notation of 3.10).

In the setup of 3.8, $\Th_c$ consists of 
$$\align&(1,S_5),(1,H_{41}),(1,H_{32}),(1,H_{311}),(1,H_{221}),(1,H_{2111}),(\ti H_{2111},H_{32}),\\&
(\ti H_{2111},H_{221}), (\ti H_{311},H_{32}), (\ti H_{221},\ti H),(H_{221},H_{221}),(H_{32},H_{32}),\\&
(H_{2111},H_{2111}), (H_{311},H_{311}),(H_{41},H_{41}),(S_5,S_5),(\ti H,\ti H),\endalign$$
(notation of 3.10).

In each case, the number of $\cg_c$-orbits on $\Th_c$ is equal to $|c|$.
By computation we see that $S_{H, H'}$ with $(H,H')$ running through a set of representatives
for the $\cg_c$-orbits on $\Th_c$ are of the form $\un{E_{H,H'}}$ (see 0.3) where 
$E_{H,H'}\in\car_c$ runs through the objects of $\car_c$ described by the rows of the matrix $M_c$
in 3.2-3.8 (in the same order as the one used in the description of $\Th_c$ given above).
These objects form a basis of $\cg_c$, due to the form of the matrix $M_c$. Now theorem 0.4 follows in our
case.

\head 4. Proof of Theorem 0.7\endhead
\subhead 4.1\endsubhead
Let $H\sub H'$ be subgroups of the finite group $\G$ with $H$ normal in $H'$.
For any $x\in\G$ we consider the set $S(x)$ of all $\mu$ in $\G/H'$ such that 
for some $\g$ in $\G/H$ contained in $\mu$ we have $x\g=\g$.
Now $Z(x)$ acts on $S(x)$ by $y:\mu\m y\mu$. For any $(x,\s)\in M(\G)$ let
$N_{x,\s}\in\NN$ be the multiplicity of $\s$ in the permutation 
representation of $Z(x)$ on $S(x)$. We have 
$$N_{x,\s}=|Z(x)|\i\sum_{y\in Z(x)}\sha(\mu\in S(x);y\mu=\mu)\tr(y,\s)$$
where 
$$\align&\sha(\mu\in S(x);y\mu=\mu)\\&
=\sha(\mu\in\G/H';\text{ for some $u\in\G$ we have }xuH=uH,\mu=uH',yuH'=uH').\endalign$$
If the previous three equations hold for some $u$ then they hold for $uh'$ for any $h'\in H'$. (Indeed,
$xuh'H=uh'H$ since $h'H=Hh'$, and $\mu=uh'H',yuh'H'=uh'H'$.) Thus,
$$\sha(\mu\in S(x);y\mu=\mu)=\sha(u\in\G;xuH=uH,yuH'=uH')/|H'|$$
and
$$\align&N_{x,\s}=|Z(x)|\i|H'|\i\sum_{y\in Z(x)}\sha(u\in\G;xuH=uH,yuH'=uH')\tr(y,\s)\\&
=|Z(x)|\i|H'|\i\sum_{y\in Z(x),u\in\G;xuH=uH,yuH'=uH'}\tr(y,\s).\endalign$$
Let $f=\sum_{(x,\s)\in M(\G)}N_{x,\s}(x,\s)\in \CC[M(\G)]$. We have $f=S_{H,H'}$.
We write $A(f)=\sum_{(x',\s')\in M(\G)}N'_{x',\s'}(x',\s')$ with $N'_{x',\s'}\in\CC$.
We have
$$\align&N'_{x',\s'}=\sum_{(x,\s)\in M(\G)}N_{x,\s}{(x,\s),(x',\s')}\\&=
\sum_{(x,\s)\in M(\G)}|Z(x)|\i|H'|\i|Z(x)|\i|Z(x')|\i
\sum_{y\in Z(x),u\in\G;xuH=uH,yuH'=uH'}\\&
\sum_{z\in\G;zxz\i x'=x'zxz\i}\ov{\tr(zxz\i,\s')}\tr(z\i x'z,\s)\tr(y,\s)  \\&
=\sum_{x\in\G}|\G|\i|H'|\i|Z(x)|\i|Z(x')|\i
\sum_{y\in Z(x),u\in\G;xuH=uH,yuH'=uH'}\\&
\sum_{z\in\G;zxz\i x'=x'zxz\i}\ov{\tr(zxz\i,\s')}
\sum_{\s\in Irr(Z(x)}\tr(z\i x'z,\s)\tr(y,\s).\endalign$$
The last sum over $\s$ equals $|Z(x)\cap Z(y)|$ if 
$z\i x'z=ay\i a\i$ for some $a\in Z(x)$ and equals $0$ otherwise.
Hence
$$\align&N'_{x',\s'}=\sum_{x\in\G}|\G|\i|H'|\i|Z(x)|\i|Z(x')|\i
\sum_{y\in Z(x),u\in\G;xuH=uH,yuH'=uH'}\\&
\sum_{z\in\G;zxz\i x'=x'zxz\i, a\in Z(x),z\i x'z=ay\i a\i}\ov{\tr(zxz\i,\s')}.\endalign$$
We substitute $z_1=za$. We get
$$\align&N'_{x',\s'}=\sum_{x\in\G}|\G|\i|H'|\i|Z(x)|\i|Z(x')|\i
\sum_{y\in Z(x),u\in\G;xuH=uH,yuH'=uH'}\\&
\sum_{z_1\in\G;z_1xz_1\i x'=x'z_1xz_1\i,a\in Z(x),z_1\i x'z_1=y\i}\ov{\tr(z_1xz_1\i,\s')}.\endalign$$
We can eliminate $a$ and change $z_1$ to $z$. We get
$$\align&N'_{x',\s'}=\sum_{x\in\G}|\G|\i|H'|\i|Z(x')|\i
\sum_{y\in Z(x),u\in\G;xuH=uH,yuH'=uH'}\\&
\sum_{z\in\G;zxz\i x'=x'zxz\i,z\i x'z=y\i}\ov{\tr(zxz\i,\s')}.\endalign$$
We substitute $x_1=u\i xu,y_1=u\i yu,z_1=zu$. We get
$$\align&N'_{x',\s'}=\sum_{x_1\in\G}|\G|\i|H'|\i|Z(x')|\i
\sum_{y_1\in Z(x_1),u\in\G;x_1H=H,y_1H'=H'}\\&
\sum_{z_1\in\G;z_1x_1z_1\i x'=x'z_1x_1z_1\i,z_1\i x'z_1=y_1\i}
\ov{\tr(z_1x_1z_1\i,\s')}.\endalign$$
We can eliminate $u$ and change $x_1,y_1,z_1$ to $x,y,z$. We get
$$\align&N'_{x',\s'}=|H'|\i|Z(x')|\i\sum_{x\in H,y\in Z(x)\cap H'}\\&
\sum_{z\in\G;zxz\i x'=x'zxz\i,z\i x'z=y\i}\ov{\tr(zxz\i,\s')}.\endalign$$
Here the condition $zxz\i x'=x'zxz\i$ follows from $z\i x'z=y\i$, $yx=xy$. Hence
$$N'_{x',\s'}=|H'|\i|Z(x')|\i\sum_{x\in H,y\in Z(x)\cap H'}
\sum_{z\in\G;z\i x'z=y\i}\ov{\tr(zxz\i,\s')},$$
that is,
$$N'_{x',\s'}
=|H'|\i|Z(x')|\i\sum_{x\in H}\sum_{z\in\G;z\i x'z\in Z(x)\cap H'}\ov{\tr(zxz\i,\s')}.$$
We substitute $zxz\i=x_1$. We get
$$N'_{x',\s'}=|H'|\i|Z(x')|\i
\sum_{x_1\in\G,z\in\G;x'\in Z(x_1)\cap zH'z\i,x_1\in zHz\i}\ov{\tr(x_1,\s')},$$
that is
$$\align&N'_{x',\s'}=|H'|\i|Z(x')|\i\\&
\sum_{z\in\G;z\i x'z\in H'}\sum_{x_1\in Z(x')\cap zHz\i}\ov{\tr(x_1,\s')}\endalign$$
and
$$\align&N'_{x',\s'}=|H'|\i|Z(x')|\i\\&
\sum_{z\in\G;z\i x'z\in H'}(1:\s'|(Z(x')\cap zHz\i))|Z(x')\cap zHz\i|\endalign$$
where $:$ denotes multiplicity. Thus we have
$$N_{x',\s'}\in\QQ_{\ge0}$$
so that $A(f)\in M(\G)_{\ge0}$. Since $f\in M(\G)_{\ge0}$ is obvious we see that $f$ is bipositive.
This proves Theorem 0.7.

\enddocument